\documentclass[12pt,oneside,reqno]{amsart}
\usepackage{}
\usepackage{amssymb}

\usepackage{bbm}
\usepackage{cases}
\usepackage{amsmath}
\usepackage{graphicx}
\usepackage{mathrsfs}
\usepackage{stmaryrd}
\usepackage{color}
\usepackage{soul}
\usepackage[dvipsnames]{xcolor}
\usepackage{amsfonts}
\usepackage{enumerate,amsmath,amssymb,amsthm}
\usepackage[colorlinks,citecolor=blue,urlcolor=blue]{hyperref}

\numberwithin{equation}{section}

\newcommand{\be}{\begin{eqnarray}}
\newcommand{\mE}{\end{eqnarray}}
\newcommand{\ce}{\begin{eqnarray*}}
\newcommand{\de}{\end{eqnarray*}}
\newtheorem{theorem}{Theorem}[section]
\newtheorem{lemma}[theorem]{Lemma}
\newtheorem{remark}[theorem]{Remark}
\newtheorem{definition}[theorem]{Definition}
\newtheorem{proposition}[theorem]{Proposition}
\newtheorem{example}[theorem]{Example}
\newtheorem{corollary}[theorem]{Corollary}

\def\eps{\varepsilon}

\def\p{\partial}

\def\[{{\Big[}}
\def\]{{\Big]}}
\def\<{{\langle}}
\def\>{{\rangle}}
\def\({{\Big(}}
\def\){{\Big)}}

\def\bx{{\mathbf{x}}}

\def\dif{{\mathord{{\rm d}}}}

\def\no{\nonumber}
\def\={&\!\!=\!\!&}
\def\bt{\begin{theorem}}
\def\et{\end{theorem}}
\def\bl{\begin{lemma}}
\def\el{\end{lemma}}
\def\br{\begin{remark}}
\def\er{\end{remark}}

\def\bd{\begin{definition}}
\def\ed{\end{definition}}
\def\bp{\begin{proposition}}
\def\ep{\end{proposition}}
\def\bc{\begin{corollary}}
\def\ec{\end{corollary}}
\def\bx{\begin{example}}
\def\ex{\end{example}}

\def\cG{{\mathcal G}}

\def\cI{{\mathcal I}}

\def\cK{{\mathcal K}}
\def\cL{{\mathcal L}}

\def\cO{{\mathcal O}}

\def\cQ{{\mathcal Q}}
\def\cR{{\mathcal R}}

\def\mE{{\mathbb E}}

\def\mN{{\mathbb N}}

\def\mR{{\mathbb R}}

\def\sL{{\mathscr L}}

\def\geq{\geqslant}
\def\leq{\leqslant}

\allowdisplaybreaks

\begin{document}

\title{Diffusion approximation for fully coupled stochastic differential equations}

\date{}

\author{}

\author{Michael R\"ockner \,\, and\,\, Longjie Xie}

\address{Michael R\"ockner:
	Fakult\"{a}t f\"{u}r Mathematik, Universit\"{a}t Bielefeld, D-33501 Bielefeld, Germany, and Academy of Mathematics and Systems Science,
	Chinese Academy of Sciences (CAS), Beijing, 100190, P.R.China\\
	Email: roeckner@math.uni-bielefeld.de
}

\address{Longjie Xie:
	School of Mathematics and Statistics and
 Research Institute of Mathematical Science, Jiangsu Normal University,
	Xuzhou, Jiangsu 221000, P.R.China\\
	Email: longjiexie@jsnu.edu.cn
}

\thanks{
This work is supported by the DFG through CRC 1283, the Alexander-von-Humboldt foundation, NSFC
(No. 11701233, 11931004) and NSF of Jiangsu (BK20170226)}

\begin{abstract}
	We consider a Poisson equation in $\mathbb R^d$ for the elliptic operator corresponding to an ergodic diffusion process. Optimal regularity and smoothness with respect to  the parameter are obtained under mild conditions on the coefficients. The result is then applied to establish a general  diffusion approximation for fully coupled multi-time-scales stochastic differential equations with only H\"older continuous coefficients. Four different averaged equations as well as rates of convergence are obtained.  Moreover, the convergence is shown to rely only on the regularities of the coefficients  with respect to  the slow variable, and does not depend on their regularities  with respect to the fast component.
	
	\bigskip

	
\end{abstract}

\maketitle

\section{Introduction}

\subsection{Poisson equation in the whole space}
The first topic of this paper is to study the following Poisson equation in $\mR^{d_1}$ ($d_1\geq1$):
\begin{align}
\sL_0(x,y)u(x,y)=f(x,y),\quad x\in\mR^{d_1},  \label{pde0}
\end{align}
where $y\in\mR^{d_2}$ ($d_2\geq 1$) is a parameter, and
\begin{align}
\sL_0:=\sL_0(x,y):=\sum_{i,j=1}^{d_1} a^{ij}(x,y)\frac{\p^2}{\p x_i\p x_j}+\sum_{i=1}^{d_1}b^i(x,y)\frac{\p}{\p x_i}.\label{l0}
\end{align}
The Poisson equation is one of the well-known equations in mathematical physics. When the above equation is formulated on a compact set, the corresponding theory is well known, see e.g. \cite{GM} or \cite[Chapter 12]{EK}. However, equation (\ref{pde0}) {\it with a parameter and in the whole space $\mR^{d_1}$} has been studied only relatively recently,  and it turns out to be
one of the key tools in
the theory of stochastic averaging, homogenization  and other limit theorems in probability theory (see e.g. \cite{FW,HP,PSV,Par}). This understanding was the reason for a series
of papers by Pardoux  and Veretennikov \cite{P-V,P-V2,P-V3}, where equation (\ref{pde0}) was first studied and then used to establish diffusion approximations for slow-fast stochastic differential equations (SDEs for short), see also \cite{HH,HMP,KV,RSX,Ve}   for further generalizations.

\vspace{2mm}
Since there is no boundary condition in equation (\ref{pde0}), the solution $u$ turns out to be defined up to an additive constant, which is quite natural due to $\sL_0(x,y)1\equiv0$. To fix this constant, it is necessary to make the following ``centering" assumption  on the potential term $f$:
\begin{align}\label{cen2}
\int_{\mR^{d_1}}f(x,y)\mu^y(\dif x)=0,\quad\forall y\in\mR^{d_2},
\end{align}
where $\mu^y(\dif x)$ is the (unique) invariant measure of an ergodic Markov process $X_t^y$ (see (\ref{sde2}) below) with a generator $\sL_0(\cdot,y)$ given $y$.
Such kind of assumption is  analogous to the centering in the standard central limit theorem (CLT for short). We mention that the main problem addressed in \cite{P-V,P-V2,P-V3} concerning equation (\ref{pde0}) is the second order regularity of the solution $u$ with respect to $x$ as well as the parameter $y$ (which is more difficult !), which would suffice for the application of It\^o formula with some diffusions plugged in for both variables,  and \cite{P-V} used mainly probability arguments while \cite{P-V2,P-V3} used  essentially results from the partial differential equations (PDEs for short) theory.
Later on, the results of \cite{P-V,P-V2,P-V3} for the Poisson equation (\ref{pde0})  have also been adopted to study the CLT, moderate and large derivations  (see e.g. \cite{CCG,DS,MS,S1}), spectral methods, averaging principle and homogenization for multi-scale systems (see e.g. \cite{APV,Br1,RSX}) as well as numerical approximation  for time-averaging estimators and the invariant measure of SDEs or stochastic partial differential equations (SPDEs for short) (see e.g. \cite{Bk,HMP,MST,PP}).

\vspace{2mm}
One of our objectives in this paper is to further study the Poisson equation (\ref{pde0}). We develop a robust method to study the regularities of the solution $u$, especially for the smoothness with respect to the parameter $y$, which leads to simplifications and extensions of the existing results. The main result in this direction is given by {\bf Theorem \ref{popde}}. Our argument is different from  all those works mentioned above, and the assumptions on the coefficients are weaker.
In fact, instead of establishing differentiablity  of the corresponding semigroup with respect to the parameter $y$ as in \cite{P-V,P-V2,P-V3}, we shall focus on the optimal regularity of the solution $u$ with respect to the $x$ variable.
Then based on a key observation of a transfer formula (see Lemma \ref{aaav} below), we use an induction argument to show that smoothness of the solution $u$ with respect to the parameter $y$ follows directly by the optimal a priori estimate with respect to the $x$ variable, which is much simpler insofar. In addition, our method has at least three more advantages: first of all, we obtain any order of differentiablity (which is important for the asymptotic expansion analysis used in \cite{KY2,KY}) as well as H\"older continuity (which will play a crucial role below for us to study diffusion approximations) of the solution with respect to the parameter $y$ under explicit conditions on the coefficients; secondly, we
provide explicit dependence on the norms of the coefficients $a,b$ as well as the potential term $f$  involved, which might be useful for numerical analysis  and has been used essentially in \cite{RSX} to study the rate of convergence in averaging principle for two-time-scales SDEs; thirdly, our argument can  also be adopted to study equation (\ref{pde0}) in Sobolev spaces with certain $L^p$ conditions on the coefficients and  the potential term, thus leading to only weak differentiablity of the solution $u$ with respect to the parameter. Such problem has been posed in \cite{Ve} and seems difficult to be  handled by the arguments used in the previous publications. For the sake of simplicity, we do not
deal with this setting in the present article and postpone it to further studies.

\subsection{Diffusion approximations}
The result of Theorem \ref{popde} will then be used to study the asymptotic problem for fully coupled multi-scale stochastic dynamical systems, which is the second topic of this paper. More precisely, consider the following in-homogeneous multi-time-scales SDE in $\mR^{d_1+d_2}$:
\begin{equation} \label{sde0}
\left\{ \begin{aligned}
&\dif X^{\eps}_t =\alpha_\eps^{-2}b(X^{\eps}_t,Y^{\eps}_t)\dif t+\beta_\eps^{-1}c(X^{\eps}_t,Y^{\eps}_t)\dif t+\alpha_\eps^{-1}\sigma(X^{\eps}_t,Y_t^\eps)\dif W^{1}_t,\\
&\dif Y^{\eps}_t =F(t,X^{\eps}_t, Y^{\eps}_t)\dif t+\gamma_\eps^{-1}H(t,X_t^\eps,Y_t^\eps)\dif t+G(t,X_t^\eps,Y_t^\eps)\dif W^{2}_t,\\
&X^{\eps}_0=x\in\mR^{d_1},\quad Y^{\eps}_0=y\in\mR^{d_2},
\end{aligned} \right.
\end{equation}
where the small parameters $\alpha_\eps,\beta_\eps,\gamma_\eps\downarrow0$ as $\eps\to0$, and without loss of generality, we may assume $\alpha_\eps^2/\beta_\eps\to0$  as $\eps\to0$.
Such  model has wide applications  in many real world dynamical systems  including planetray motion, climate
models (see e.g. \cite{K,MTV}), geophysical fluid flows (see e.g. \cite{GD}), intracellular biochemical reactions (see e.g. \cite{BKRP}), etc. We refer the interested readers to the books \cite{Ku,PS} for a more comprehensive overview. Note that we
are considering (\ref{sde0}) in the whole space and not just on compact sets. Moreover,
there exist two time-scales in the fast component $X_t^\eps$ and even the slow motion $Y_t^\eps$ has a fast
varying term. Usually, the underlying  system (\ref{sde0}) is difficult to deal with due to the widely separated time-scales and the cross interactions of slow and fast modes. Hence a simplified equation which governs the evolution of the system for small $\eps$ is highly desirable.   We also mention  that the infinitesimal generator corresponding to $(X_t^\eps, Y_t^\eps)$ has the form
$$
\sL_\eps:=\alpha_\eps^{-2}\sL_0(x,y)+\beta_\eps^{-1}\sL_3(t,x,y)+\gamma_\eps^{-1}\sL_2(t,x,y)+\sL_1(t,x,y),
$$
where $\sL_0(x,y)$ is given by (\ref{l0}) with $a(x,y):=\sigma\sigma^*(x,y)/2$, and
\begin{equation} \label{lll}
\begin{aligned}
&\sL_3:=\sL_3(x,y):=\sum_{i=1}^{d_1}c^i(x,y)\frac{\p}{\p x_i},\\
&\sL_2:=\sL_2(t,x,y):=\sum_{i=1}^{d_2}H^i(t,x,y)\frac{\p}{\p y_i},\\
&\sL_1:=\sL_1(t,x,y):=\sum_{i,j}^{d_2} \cG^{ij}(t,x,y)\frac{\p^2}{\p y_i\p y_j}+\sum_{i=1}^{d_2}F^i(t,x,y)\frac{\p}{\p y_i}\\
\end{aligned}
\end{equation}
with $\cG(t,x,y)=GG^*(t,x,y)/2$.
Thus the asymptotic behavior of SDE (\ref{sde0})  as $\eps\to0$  is closely related to the limit theorem for solutions of second order parabolic and elliptic equations with singularly perturbed terms, which has its own interest in the theory of PDEs, see e.g. \cite{HP,Par} and \cite[Chapter IV]{Fr}.

\vspace{2mm}
When $c=H\equiv0$, the celebrated theory of averaging principle asserts that the slow motion $Y_t^\eps$ will convergence  in distribution as $\eps\to 0$ to the solution $\bar Y_t$ of the following reduced equation in $\mR^{d_2}$:
\begin{align}\label{sde00}
\dif \bar Y_t=\bar F(t, \bar Y_t)\dif t+\bar G(t,\bar Y_t)\dif W^2_t,\quad \bar Y_0=y,
\end{align}
where the new averaged coefficients are given by
$$
\bar F(t,y):=\!\int_{\mR^{d_1}}\!F(t,x,y)\mu^y(\dif x)
$$
and
$$
 G(t,y):=\!\sqrt{\int_{\mR^{d_1}}\!G(t,x,y)G(t,x,y)^{*}\mu^y(\dif x)},
$$
and for each $y\in\mR^{d_2}$, $\mu^y(\dif x)$ is the unique invariant measure for $X_t^y$ which satisfies  the frozen equation
\begin{align}\label{sde2}
\dif X_t^y=b(X_t^y,y)\dif t+\sigma(X_t^y,y)\dif W_t^1,\quad X_0^y=x\in\mR^{d_1}.
\end{align}
The corresponding results are also known as the averaging principle of functional law of large numbers (LLN for short) type and have been  intensively
studied by the classical time discretisation method, see e.g. \cite{AV,BK,CFKM,GR,Ki,Li,Pu,V0}, see also \cite{Ce,CF,CL,WR} for similar results for SPDEs. Generalization to the general case that $\alpha_\eps=\beta_\eps=\gamma_\eps$ was first carried out by Papanicolaou, Stroock and Varadhan \cite{PSV} for a compact state space for the fast component and  time-independent coefficients, see also \cite{Ba} for a similar result in terms of PDEs.  It was found  that  the limit
distribution of the slow component will be obtained
in terms of the solution of an auxiliary Poisson equation. Such result can be regarded as an averaging principle of functional CLT type and is often  called diffusion approximation, which is  important for applications in homogenization. Later a  non-compact and homogeneous case with $c\equiv0$ and $\alpha_\eps=\gamma_\eps$ was studied in  \cite{P-V,P-V2,P-V3} by using the method of martingale problem and in \cite{KY} by the asymptotic expansion approach, see also \cite[Chapter 11]{PS}. We also mention that for numerical purposes, the existence
of the effective system (\ref{sde00}) is not enough, and the rate of convergence of the slow variable to its limit distribution has to be derived. The main motivation comes from the well-known Heterogeneous Multi-scale Methods  used to approximate the slow component $Y_t^\eps$, see e.g. \cite{Br3,Br1,EL,KVa}. However, getting error bounds is significantly harder than just showing convergence, and many of the methods commonly
employed to show distributional convergence only possibly yield a convergence rate after serious added effort.
In this direction, there are many works devoted to study the convergence rate in the averaging principle of LLN type for SDE (\ref{sde0}) when $c=H\equiv0$, see e.g. \cite{G, L1,V,ZFWL} and the references therein. As far as we know, there is still no result concerning  error bounds for the CLT type convergence of system (\ref{sde0}) in the general case.

\vspace{2mm}
We will study the asymptotic problem for system (\ref{sde0}) more systematically. The main result is given by {\bf Theorem \ref{main2}}.  It turns out that, depending on the orders how $\alpha_\eps,\beta_\eps,\gamma_\eps$ go to zero, we shall have four different regimes of interactions, which lead to four different asymptotic behaviors of system (\ref{sde0}) as $\eps\to0$, i.e.,
\begin{equation}\label{regime}
\left\{\begin{aligned} &\lim_{\eps\to0}\frac{\alpha_\eps}{\gamma_\eps}=0\quad\text{and}\quad\lim_{\eps\to0}\frac{\alpha_\eps^2}{\beta_\eps\gamma_\eps}=0,\qquad \text{Regime 1};\\
&\lim_{\eps\to0}\frac{\alpha_\eps}{\gamma_\eps}=0\quad\text{and}\quad\alpha_\eps^2=\beta_\eps\gamma_\eps,\qquad\quad\,\, \text{Regime 2};\\
&\alpha_\eps=\gamma_\eps\quad\,\,\,\quad\text{and}\quad\lim_{\eps\to0}\frac{\alpha_\eps}{\beta_\eps}=0,\qquad\quad\! \text{Regime 3};\\
&\alpha_\eps=\beta_\eps=\gamma_\eps,\qquad\qquad\qquad\quad\quad\qquad\,\,\, \text{Regime 4}.
\end{aligned}\right.
\end{equation}
If $\alpha_\eps$ and $\alpha_\eps^2$ go to zero faster than $\gamma_\eps$ and $\beta_\eps\gamma_\eps$ respectively (Regime 1), we show that the limit behavior of system (\ref{sde0}) coincides with SDE (\ref{sde00}), i.e., the traditional case  corresponding to $c=H\equiv0$; if $\alpha_\eps$ goes to zero faster than $\gamma_\eps$ while $\alpha_\eps^2$ and $\beta_\eps\gamma_\eps$ are of the same order (Regime 2),
then the homogenization effect of term $c$ will occur  in the limit dynamics; whereas if $\alpha_\eps$ and $\gamma_\eps$ are of the same order and $\alpha_\eps$ goes
to zero faster than $\beta_\eps$ (Regime 3), then the homogenization effect of term $H$ appears; finally, when all the parameters are of the same order (Regime 4), then homogenization effects of term $c$ and term $H$ will occur together.

\vspace{2mm}
We shall handle Regime 1-4 in a robust and unified way. Our method relies only on the technique of Poisson equation (\ref{pde0}), and does not involve extra time discretisation procedure (see e.g. \cite{BK,CL,Li,V0}), martingale problem (see \cite{MS,PSV,P-V,P-V2,P-V3,S1}) nor asymptotic expansion argument (see \cite{KY2,KY}) and thus is quite simple. Moreover, the conditions on the coefficients are  weaker (only H\"older continuous) than the known results in the literature, and rates of convergence are obtained as easy by-products of our argument,  which we believe are rather sharp, see  {\bf Remark \ref{br1}} and {\bf Remark \ref{br2}} for more explanations. We also point out that unlike the above mentioned results, where the second order regularity of the solution to the Poisson equation with respect to the parameter is commonly needed in the proof, we only use its H\"older continuity, which also simplifies the arguments used in \cite{RSX} and appear intuitively natural, since the second order derivative of the solution to the Poisson equation with respect to the parameter does not appear in the final limit equation.  Throughout our proof,
two new fluctuation estimates of functional LLN type in Lemma \ref{key} and functional CLT type in Lemma \ref{key22} will play important roles, which might be used to study other limit theorems and should be of independent interest.

\vspace{2mm}
The rest of this paper is organized as follows. In Section 2 we provide the assumptions and state our main results. Section 3 is devoted to the study of Possion equation (\ref{pde0}) and we \ prove Theorem \ref{popde}. In Section 4 we prepare two fluctuation lemmas, and then we give the proof of Theorem \ref{main2} in Section 5. Throughout our paper, we use the following convention: $C$ and $c$ with or without subscripts will denote positive constants, whose values may change in different places, and whose dependence on parameters can be traced from the calculations.

\vspace{2mm}
{\bf Notations:} To end this section, we introduce some notations. Let $\mN^*:=\{1,2\cdots\}$. Given a function space, the subscript $b$ will stand for boundness, while the subscript $p$ stands for polynomial growth in $x$. More precisely, for a function $f(t,x,y)\in L^\infty_p:=L^{\infty}_p(\mR_+\times\mR^{d_1+d_2})$, we mean there exist constants $C, m>0$ such that
$$
|f(t,x,y)|\leq C(1+|x|^m),\quad\forall t>0, x\in\mR^{d_1}, y\in\mR^{d_2}.
$$
For $0<\delta\leq 1$, the space $C_p^{\delta,0}:=C_p^{\delta,0}(\mR^{d_1+d_2})$ consists of all functions which are local H\"older continuous and have at most polynomial growth  in $x$ uniformly with respect to $y$, i.e.,  there exist constants $C, m>0$ such that for any $x_1,x_2\in\mR^{d_1}$,
$$
|f(x_1,y)-f(x_2,y)|\leq C\big(|x_1-x_2|^\delta\wedge1\big)\big(1+|x_1|^m+|x_2|^m\big),\quad\forall y\in\mR^{d_2}.
$$
We also define a quasi-norm for $C^{\delta,0}_p$ by
$$
[f]_{C^{\delta,0}_p}:=\sup_{y\in\mR^{d_2}}\sup_{|x_1|,|x_2|\leq 1,|x_1-x_2|\leq 1}|f(x_1,y)-f(x_2,y)|/|x_1-x_2|^\delta.
$$
For $0<\vartheta< 1$, the space $C_p^{\delta,\vartheta}:=C_p^{\delta,\vartheta}(\mR^{d_1+d_2})$ consists of all functions that are $\delta$-local H\"older continuous with polynomial growth in $x$ and  $\vartheta$-H\"older continuous in $y$, i.e.,  there exist constants $C, m>0$ such that for any $x_1,x_2\in\mR^{d_1}$ and $y_1,y_2\in\mR^{d_2}$,
\begin{align*}
|f(x_1,y_1)-f(x_2,y_2)|&\leq C\Big[\big(|x_1-x_2|^\delta\wedge1\big)+\big(|y_1-y_2|^\vartheta\wedge1\big)\Big]\\
&\qquad\times\big(1+|x_1|^m+|x_2|^m\big).
\end{align*}
Similarly, we define a quasi-norm for $C^{\delta,\vartheta}_p$ by
$$
[f]_{C^{\delta,\vartheta}_p}:=\sup_{|y_1-y_2|\leq 1}\sup_{|x_1|,|x_2|\leq 1,|x_1-x_2|\leq 1}\frac{|f(x_1,y_1)-f(x_2,y_2)|}{|x_1-x_2|^\delta+|y_1-y_2|^\vartheta}.
$$
When $\delta,\gamma\geq1$, we use $C_p^{\delta,\vartheta}:=C_p^{\delta,\vartheta}(\mR^{d_1+d_2})$ to denote the space of all functions $f$ satisfying $\p^{[\delta]}_x\p^{[\vartheta]}_yf\in C_p^{\delta-[\delta],\vartheta-[\vartheta]}$. Finally, $ C_p^{\gamma,\delta,\vartheta}:=C_p^{\gamma,\delta,\vartheta}(\mR_+\times\mR^{d_1+d_2})$ with $0<\gamma\leq 1$ denotes the space of all functions $f$ such that  for every fixed $t>0$, $f(t,\cdot,\cdot)\in C_p^{\delta,\vartheta}$ and for every $(x,y)\in \mR^{d_1+d_2}$, $f(,\cdot,x,y)\in C_b^\gamma(\mR_+)$, where $C_b^\gamma$ is the usual bounded H\"older space.

\section{Assumptions and main results}

To state our main results, we first introduce some basic assumptions. Throughout this paper, we shall always assume the following non-degeneracy conditions on the diffusion coefficients:

\vspace{2mm}
\noindent{\bf (A$_\sigma$):} the coefficient $a=\sigma\sigma^*$ is non-degenerate in $x$ uniformly with respect to $y$, i.e., there exists $\lambda>1$ such that  for any $y\in\mR^{d_2}$,
$$
\lambda^{-1}|\xi|^2\leq |a(x,y)\xi|^2\leq \lambda|\xi|^2,\ \ \forall\xi\in\mR^{d_1}.
$$

\vspace{1mm}
\noindent{\bf (A$_{G}$):} the coefficient $\cG=GG^*$ is non-degenerate in $y$ uniformly with respect to $(t,x)$, i.e., there exists $\lambda>1$ such that for any $(t,x)\in\mR_+\times\mR^{d_1}$,
$$
\lambda^{-1}|\xi|^2\leq |\cG(t,x,y)\xi|^2\leq\lambda|\xi|^2,\ \ \forall\xi\in\mR^{d_2}.
$$

\vspace{1mm}
\noindent  Note that the operator $\sL_0$ in the Poisson equation (\ref{pde0}) can be viewed as the infinitesimal generator of the frozen SDE (\ref{sde2}). We make the following very weak recurrence assumption on the drift $b$  to ensure the existence of an  invariant measure $\mu^y(\dif x)$ for $X_t^y$:

\vspace{3mm}
\noindent{\bf (A$_b$):}\qquad\qquad\qquad\quad $\lim_{|x|\to\infty}\sup_y \<x,b(x,y)\>=-\infty$.

\vspace{2mm}
Our main result concerning the Poisson equation (\ref{pde0}) is as follows.

\bt\label{popde}
Let {\bf (A$_\sigma$)} and {\bf (A$_b$)} hold. Assume that $a, b\in C_b^{\delta,\eta}$ with $0<\delta\leq 1$ and $\eta\geq0$. Then for every function $f\in C_p^{\delta,\eta}$ satisfying (\ref{cen2}), there exists a unique solution $u\in C_p^{2+\delta,\eta}$ to equation (\ref{pde0}) which also satisfies (\ref{cen2}).
Moreover, there exist constants $m>0$ and $C_0>0$ depending only on $d_1,d_2$ and $\|a\|_{C_b^{\delta,0}}, \|b\|_{C_b^{\delta,0}}, [f]_{C_p^{\delta,0}}$ such that:\\
(i) (Case $\eta=0$) for any $x\in\mR^{d_1}$ and $y\in\mR^{d_2}$,
\begin{align}\label{key1}
|u(x,y)|+|\nabla_xu(x,y)|+|\nabla_x^2u(x,y)|\leq C_0(1+|x|^m),
\end{align}
and for any $x_1,x_2\in\mR^{d_1}$,
\begin{align}\label{key11}
|\nabla_x^2u(x_1,y)-\nabla_x^2u(x_2,y)|\leq C_0\big(|x_1-x_2|^\delta\wedge1\big)(1+|x_1|^m+|x_2|^m);
\end{align}
(ii) (Case $\eta>0$) for any $x\in\mR^{d_1}$,
\begin{align}\label{key3}
\|u(x,\cdot)\|_{C_b^\eta}\leq C_0\cK_\eta(1+|x|^m),
\end{align}
where  $\cK_\eta>0$ is a constant depending only on $\|a\|_{C_b^{\delta,\eta}}, \|b\|_{C_b^{\delta,\eta}}$ and $[f]_{C_p^{\delta,\eta}}$, which is defined recursively by (\ref{dep}) below.
\et

\br
(i) Estimates (\ref{key1}) and (\ref{key11}) reflect the optimal regularity of the solution $u$ with respect to the $x$ variable, which are natural because we can get the two derivatives  for free by virtue of the elliptic property of the operator. Note that we allow the potential term $f$ to have polynomial growth in $x$. Estimate (\ref{key3}) means that we need both the coefficients and the right hand side $f$ to be  differentiable with respect to the $y$ variable in order to guarantee the same regularity for $u$ with respect to $y$, which is also reasonable and optimal  since $y$ is only a parameter in the equation. Note that only H\"older continuity of the coefficients with respect to the $x$ variable is needed. The key point of our arguments which differs from \cite{P-V,P-V2,P-V3} is that we show that the smoothness of the solution with respect to the parameter $y$ follows directly from its optimal regularity with respect to the $x$ variable.

\vspace{1mm}
\noindent(ii) The positive constant $\cK_\eta$ in (\ref{key3}) can be given explicitly in terms of $\|a\|_{C_b^{\delta,\eta}}, \|b\|_{C_b^{\delta,\eta}}$ and $[f]_{C_p^{\delta,\eta}}$. In fact, let
$$
\kappa_\eta:=\|a\|_{C_b^{\delta,\eta}}+\|b\|_{C_b^{\delta,\eta}}\quad\text{and}\quad\cK_0:=1.
$$
Then when $\eta$ is an integer, we have
\begin{align}\label{dep}
\cK_\eta:=[f]_{C_p^{\delta,\eta}}+\sum_{\ell=1}^{\eta}C_\eta^\ell\cdot \kappa_\ell\cdot\cK_{\eta-\ell},
\end{align}
while for $\eta\in(0,\infty\setminus\mN^*)$, we have
\begin{align}\label{dep2}
\cK_\eta=\|f\|_{C^{\delta,\eta}_p}+\sum_{\ell=1}^{[\eta]} C_{[\eta]}^\ell\cdot\kappa_\ell \cdot\cK_{\eta-\ell}+\sum_{\ell=0}^{[\eta]}C_{[\eta]}^{[\eta]-\ell}\cdot\kappa_{\eta-\ell}\cdot\cK_{\ell}.
\end{align}
This will be useful for numerical purposes to approximate the solution of the Poisson equation with singular coefficients, see e.g. \cite{RSX}.

\vspace{1mm}
\noindent(iii) The Poisson equation considered in \cite{P-V,P-V3} is not fully coupled, i.e., the operator $\sL_0$ does not depends on the parameter $y$. Growth conditions are imposed on the term $f$ and the Sobolev regularity of the solution $u$ with respect to the $x$ variable is obtained. While, \cite{P-V3} mainly focus on the degenerate case. In \cite{P-V2}, the Poisson equation (\ref{pde0}) was investigated
in  H\"older classes of functions as in the present paper. Here, we relax the regularity assumptions on the  coefficients   and  improve the regularity results for the solution.
\er

Now we turn to the multi-time-scales stochastic dynamical system (\ref{sde0}).
As shown in \cite{PSV}, the limit behavior for the slow component $Y_t^\eps$ will be given in terms of the solution of an auxiliary Poisson equation involving the drift $H$. Thus, it is necessary to make the following assumption:

\vspace{1mm}
\noindent{\bf (A$_H$):} the drift $H$ is centered, i.e.,
\begin{align}\label{cen}
\int_{\mR^{d_1}}H(t,x,y)\mu^y(\dif x)=0,\quad\forall (t,y)\in\mR_+\times\mR^{d_2},
\end{align}
where $\mu^y(\dif x)$ is the invariant measure for SDE (\ref{sde2}).

\vspace{1mm}
Note that the drift $c$ is not involved in the frozen equation (\ref{sde2}). We need the following additional condition for $c$ to ensure the non-explosion of the solution $X_t^\eps$: for $\eps>0$ small enough, it holds that
\begin{align}\label{ac}
\lim_{|x|\to\infty}\sup_y \<x,b(x,y)+\eps c(x,y)\>=-\infty.
\end{align}

Under (\ref{cen}) and according to Theorem \ref{popde}, there exists a unique solution $\Phi(t,x,y)$ to the following Poisson equation in $\mR^{d_1}$:
\begin{align}
\sL_0(x,y)\Phi(t,x,y)=-H(t,x,y),\quad x\in\mR^{d_1},  \label{pde1}
\end{align}
where $(t,y)\in\mR_+\times\mR^{d_2}$ are regarded as parameters.
We introduce the new averaged drift coefficients by
\begin{align*}
\hat F_1(t,y)&:=\int_{\mR^{d_1}}\!F(t,x,y)\mu^{y}(\dif x);\\
\hat F_2(t,y)&:=\int_{\mR^{d_1}}\!\big[F(t,x,y)+c(x,y)\cdot\nabla_x\Phi(t,x,y)\big]\mu^{y}(\dif x);\\
\hat F_3(t,y)&:=\int_{\mR^{d_1}}\!\big[F(t,x,y)+H(t,x,y)\cdot\nabla_y\Phi(t,x,y)\big]\mu^{y}(\dif x);\\
\hat F_4(t,y)&:=\int_{\mR^{d_1}}\!\big[F(t,x,y)+c(x,y)\cdot\nabla_x\Phi(t,x,y)\\
&\qquad\qquad\qquad\qquad\!+H(t,x,y)\cdot\nabla_y\Phi(t,x,y)\big]\mu^{y}(\dif x),
\end{align*}
which correspond to Regime 1-Regime 4 described in (\ref{regime}), and the new diffusion coefficients are given by
\begin{align*}
\hat G_1(t,y)&=\hat G_2(t,y):=\sqrt{\int_{\mR^{d_1}}GG^*(t,x,y)\mu^y(\dif x)};\\
\hat G_3(t,y)&=\hat G_4(t,y):=\sqrt{\int_{\mR^{d_1}}\Big[GG^*(t,x,y)+H(t,x,y)\Phi^*(t,x,y)\Big]\mu^y(\dif x)}.
\end{align*}
Hence, the precise formulation of the limit equation for SDE (\ref{sde0}) will be of the following form: for $k=1,\cdots,4$,
\begin{align}\label{sde000}
\dif \hat Y^k_t=\hat F_k(t, \hat Y^k_t)\dif t+\hat G_k(t,\hat Y^k_t)\dif W^2_t,\quad \hat Y^k_0=y.
\end{align}

The second main result of this paper is as follows.

\bt\label{main2}
Let {\bf (A$_\sigma$)}-{\bf (A$_b$)}-{\bf (A$_{G}$)}-{\bf (A$_H$)}-(\ref{ac}) hold,  $T>0$ and $\delta\in(0,1]$.

\vspace{1mm}
\noindent(i) (Regime 1) If $b, \sigma\in C_b^{\delta,\vartheta}$, $F, H, G\in C_p^{\vartheta/2, \delta,\vartheta}$ with $\vartheta\in (0,2]$, $c\in L^\infty_p$, and assume further that $\lim_{\eps\to0}\alpha_\eps^\vartheta/\gamma_\eps=0$, then for every $\varphi\in C_b^{2+\vartheta}(\mR^{d_2})$, we have
\begin{align*}
\sup_{t\in[0,T]}\Big|\mE[\varphi(Y_t^\eps)]-\mE[\varphi(\hat Y^1_t)]\Big|\leq C_T\,\Big(\frac{\alpha_\eps^{\vartheta}}{\gamma_\eps}+\frac{\alpha_\eps^2}{\gamma_\eps^2}+\frac{\alpha_\eps^2}{\beta_\eps\gamma_\eps}\Big);
\end{align*}
(ii) (Regime 2) if $b,\sigma \in C_b^{\delta,\vartheta}$,  $F, H, G\in C_p^{\vartheta/2, \delta,\vartheta}$ and $c\in C_{p}^{\delta,\vartheta}$ with $\vartheta\in (0,2]$, and assume further that $\lim_{\eps\to0}\alpha_\eps^\vartheta/\gamma_\eps=0$, then for every $\varphi\in C_b^{2+\vartheta}(\mR^{d_2})$, we have
\begin{align*}
\sup_{t\in[0,T]}\Big|\mE[\varphi(Y_t^\eps)]-\mE[\varphi(\hat Y^2_t)]\Big|\leq C_T\,\Big(\frac{\alpha_\eps^{\vartheta}}{\gamma_\eps}+\frac{\alpha_\eps^2}{\gamma_\eps^2}+\frac{\alpha_\eps^2}{\beta_\eps}\Big);
\end{align*}
(iii) (Regime 3) if $b, \sigma\in C_b^{\delta,1+\vartheta}$, $F, G\in C_{p}^{\vartheta/2,\delta,\vartheta}$, $H\in C_p^{(1+\vartheta)/2, \delta,1+\vartheta}$ with $\vartheta\in (0,1]$ and $c\in L^\infty_p$, then  for every $\varphi\in C_b^{2+\vartheta}(\mR^{d_2})$, we have
\begin{align*}
\sup_{t\in[0,T]}\Big|\mE[\varphi(Y_t^\eps)]-\mE[\varphi(\hat Y^3_t)]\Big|\leq C_T\,\Big({\alpha_\eps^\vartheta}+\frac{\alpha_\eps}{\beta_\eps}\Big);
\end{align*}
(iv) (Regime 4) if $b, \sigma\in C_b^{\delta,1+\vartheta}$, $F, G\in C_{p}^{\vartheta/2,\delta,\vartheta}$, $H\in C_p^{(1+\vartheta)/2, \delta,1+\vartheta}$ and $c\in C_{p}^{\delta,\vartheta}$ with $\vartheta\in (0,1]$,  then for every $\varphi\in C_b^{2+\vartheta}(\mR^{d_2})$, we have
\begin{align*}
\sup_{t\in[0,T]}\Big|\mE[\varphi(Y_t^\eps)]-\mE[\varphi(\hat Y^4_t)]\Big|\leq C_T\,\alpha_\eps^{\vartheta},
\end{align*}
where for $k=1,\cdots,4$, $\hat Y_t^k$ are the unique weak solutions for SDE (\ref{sde000}), and $C_T>0$ is a constant independent of $\delta, \eps$.
\et

\br\label{br1}
(i) Note that in each case, the convergence rates do not depend on the index $\delta$. This suggests that the convergence in the averaging principle for system (\ref{sde0}) relies only on the regularity of the coefficients  with respect to the $y$ (slow) variable, and does not depend on their regularity  with respect to the $x$ (fast) variable.

\vspace{1mm}

\noindent(ii)  In \cite{S1}, the weak convergence of the deviations of $Y_t^\eps$ around its averaged motion $\bar Y_t$, i.e., $Z_t^\eps:=(Y_t^\eps-\bar Y_t)/\lambda_\eps$ with $\lambda_\eps\to0$ as $\eps\to0$, was studied in a particular homogeneous case where $\alpha_\eps=\delta/\sqrt{\eps}$, $\beta_\eps=\delta$,  $\gamma_\eps=\delta/\eps$ (which is a particular case of Regime 2) and with small noise perturbations, i.e., with $G$ replaced by $\sqrt{\eps}G$  in SDE (1.4). Regime specific analysis is also done therein due to the flexibility of $\lambda_\eps$. In a very recent paper \cite{RX}, we further consider   the strong convergence and the central limit theorem for SDE (1.4) when $G$ does not depend on the $x$ variable as well as in the small noise perturbation case, which totally cover the results in \cite{S1}. The arguments developed in this paper and the Poisson equation result Theorem \ref{popde}   play important role in \cite{RX}.
\er

Let us give more explanations on the convergence rates and the assumptions on the coefficients in the above result, which we think are rather sharp.

\br\label{br2}
(i) Note that when $c=H\equiv0$, we can take $\gamma_\eps=\beta_\eps\equiv1$. In this particular case, the above result of Regime 1 simplifies to: if $b, \sigma\in C_b^{\delta,\vartheta}$ and $F, G\in C_p^{\vartheta/2, \delta,\vartheta}$ with $\vartheta\in (0,2]$, then for every $\varphi\in C_b^{2+\vartheta}(\mR^{d_2})$, we have
\begin{align*}
\sup_{t\in[0,T]}\Big|\mE[\varphi(Y_t^\eps)]-\mE[\varphi(\hat Y^1_t)]\Big|\leq C_T\,\alpha_\eps^{\vartheta}.
\end{align*}
This is known to be optimal when $\vartheta=2$. In the general case, the rate of convergence will be dominated by the fast term $\gamma_\eps^{-1}\int_0^tH(s,X_s^\eps,Y_s^\eps)\dif s$ in the slow equation.
Thus, the extra assumption that $\lim_{\eps\to0}\alpha_\eps^\vartheta/\gamma_\eps$ $ =0$ in Regime 1-2 with only H\"older coefficients is  necessary and  will be automatically satisfied when $\vartheta\geq 1$. It is also interesting to note that we may have $\lim_{\eps\to0}\alpha_\eps/\beta_\eps=+\infty$ in these two regimes.

\vspace{2mm}
(ii) It will be clear from our proof that the convergence rates of $Y_t^\eps$ to $\hat Y_t^k$ ($k=1,\cdots, 4$) depend only on the regularities of the averaged coefficients in the limit equations. Thus the regularities of the coefficients in the original equation with respect to the fast variable do not play any role, which appears to be intuitively natural. Our assumptions seem to be the weakest in order to get the desired result. Let us explain  this with Regime 4 when the coefficients are time-independent.
In order to get the $\alpha_\eps^\vartheta$-order convergence of $Y_t^\eps$ to $\hat Y_t^4$, we shall need $\hat F_4, \hat G_4\in C_b^{\vartheta}$.  Thanks to the assumption  that $b, \sigma\in C_b^{\delta,1+\vartheta}$, $H\in C_{p}^{\delta,1+\vartheta}$ and Theorem \ref{popde}, we can get that
$\Phi\in C_{p}^{2+\delta,1+\vartheta}$, which in turn means that $\nabla_y\Phi\in C_{p}^{\delta,\vartheta}$. This together with the assumptions $F,c\in C_p^{\delta,\vartheta}$ and Lemma {\ref{aaav}} yields that $\hat F_4\in C_b^\vartheta$. Similar computations hold for $\hat G_4$ as well as Regime 1-Regime 3.
\er

\section{Poisson equation in the whole space}

This section is devoted to study the Poisson equation (\ref{pde0}) in the whole space. Note that formally,
the solution $u$ should have the
following probability representation:
\begin{align}\label{pro}
u(x,y)=\int_0^\infty\mE f\big(X_t^y(x),y\big)\dif t,
\end{align}
where $X_t^y(x)$ is the unique  solution for the frozen SDE (\ref{sde2}). The assertion that $u$ satisfies (\ref{cen2}) was proven by \cite{P-V}. Let $p_t(x,x';y)$ be the density function of $X_t^y(x)$ (which is also the unique fundamental solution for the operator $\sL_0$).
For simplify, we denote by $T_tf(x,y)$ the semigroup corresponding to $X_t^y(x)$, i.e.,
$$
T_tf(x,y):=\mE \big(f(X_t^y(x),y\big)=\int_{\mR^{d_1}}p_t(x,x';y)f(x',y)\dif x'.
$$
Unlike the previous publications, we do not focus on the differentiablity of the semigroup $T_tf$ with respect to the parameter $y$. Instead, we shall first prove the optimal regularity estimates (\ref{key1}) and (\ref{key11}) for the solution $u$ with respect to the $x$-variable, then we use an induction argument to show that  estimate (\ref{key3}) with respect to the parameter $y$ follows directly.  Throughout this section, we shall always assume that
{\bf (A$_\sigma$)}-{\bf (A$_b$)} hold.

In order to study the optimal regularity for equation (\ref{pde0}), let us first collect some classical results concerning the properties for $p_t(x,x';y)$.

\bl
Assume {\bf (A$_\sigma$)} holds and $T>0$. Let $a,b\in C_b^{\delta,0}$ with  $0<\delta\leq 1$. Then
for every $\ell=0,1,2$ and any $0<t\leq T$, we have
\begin{align}\label{ppp1}
|\nabla_x^\ell p_t(x,x';y)|\leq C_{T}t^{-(d+\ell)/2}\exp\big(-c_0|x-x'|^2/t\big),
\end{align}
and
\begin{align}\label{ppp3}
\bigg|\int_{\mR^{d_1}}\nabla_x^2 p_t(x,x';y)\dif x'\bigg|\leq C_{T}t^{-(2-\delta)/2},
\end{align}
and for every $x_1,x_2\in\mR^{d_1}$,
\begin{align}\label{ppp2}
\begin{split}
&|\nabla_x^2 p_t(x_1,x';y)-\nabla_x^2 p_t(x_2,x';y)|\\
&\leq C_{T}\Big[\big(|x_1-x_2|\wedge1\big)t^{-(d+3)/2}+\big(|x_1-x_2|^{\delta}\wedge1\big)t^{-(d+2)/2}\Big]\\
&\quad\times\Big(\exp\big(\!-\!c_0|x_1-x'|^2/t\big)+\exp\big(\!-\!c_0|x_2-x'|^2/t\big)\Big),
\end{split}
\end{align}
where $C_T,c_0>0$ are constants independent of $y$.

If we further assume {\bf (A$_b$)} holds, then  the limit
$$
p_\infty(x',y):=\lim_{t\to\infty}p_t(x,x';y)
$$
exists, and for every $k, j\in\mR_+$, there exists a constant $m>0$ such that for any $t\geq 1$, $x, x'\in\mR^{d_1}$ and $y\in\mR^{d_2}$,
\begin{align}\label{pin1}
|p_\infty(x',y)|\leq \frac{C_0}{1+|x'|^j},
\end{align}
and
\begin{align}\label{pin2}
|p_t(x,x';y)-p_\infty(x',y)|\leq C_0\frac{1+|x|^m}{(1+t)^k(1+|x'|^j)}.
\end{align}
where  $C_0$ is a positive constant depending only on $\lambda,d_1,d_2$ and $\|a\|_{C_b^{\delta,0}}, \|b\|_{C_b^{\delta,0}}$.
\el
\begin{proof}
	Estimate (\ref{ppp1}) is well-known, see e.g. \cite[Chapter IV, $\S$13, (13.1)]{La-So-Ur}, while estimate (\ref{ppp3}) can be found in \cite[Chapter IV, $\S$14, (14.2)]{La-So-Ur}. When $|x_1-x_2|\leq 1$, estimate (\ref{ppp2}) follows by \cite[Chapter IV, $\S$13, (13.2)]{La-So-Ur}. While for $|x_1-x_2|>1$, the conclusion follows easily by (\ref{ppp1}).
	Finally, estimates (\ref{pin1}) and (\ref{pin2}) were given in \cite[Proposition 3]{P-V2}.
\end{proof}

To shorten the notation, we will write for $\ell\in\mN^*$,
$$
\frac{\p^\ell\sL_0}{\p y^\ell}(x,y):=\sum_{i,j=1}^{d_1} \p_y^\ell a^{ij}(x,y)\frac{\p^2}{\p x_i\p x_j}+\sum_{i=1}^{d_1}\p_y^\ell b^i(x,y)\frac{\p}{\p x_i}.
$$
Given a function $h(x,y)$ on $\mR^{d_1}\times\mR^{d_2}$, we shall denote by $\bar h(y)$ its average with respect to the measure $\mu^y(\dif x)$, i.e.,
\begin{align}\label{aver}
\bar h(y):=\int_{\mR^{d_1}}h(x,y)\mu^y(\dif x).
\end{align}
Then it is easy to see that
\begin{align}\label{add}
f(x,y):=h(x,y)-\bar h(y)
\end{align}
satisfies the centering condition (\ref{cen2}).	The following result will play an important role in the study of the smoothness of the solution to the Poisson equation with respect to the parameter $y$ as well as the diffusion approximations for SDE (\ref{sde0}).

\bl\label{aaav}
Let {\bf (A$_\sigma$)} and {\bf (A$_b$)} hold, $a, b\in C_b^{\delta,\eta}$ with $0<\delta\leq 1$ and $\eta> 0$. Given a function $h\in C_p^{\delta,\eta}$, let  $f$ be defined by (\ref{add}). Assume that there exists a solution $u\in C_p^{2+\delta,(\eta-1)\vee0}$ to the Poison equation (\ref{pde0}). Then,

\vspace{1mm}
\noindent(i) if $\eta\in\mN^*$,  we have
\begin{align}\label{trans}
\p^{\eta}_y\bar h(y)=\int_{\mR^{d_1}}\!\Big[\p_y^{\eta} h(x,y)-\sum_{\ell=1}^{\eta} C_{\eta}^\ell \cdot\frac{\p^\ell\sL_0}{\p y^\ell}(x,y)\p_y^{\eta-\ell}u(x,y)\Big]\mu^y(\dif x);
\end{align}
(ii) if $\eta\in(0,\infty)\setminus\mN^*$, we have that for any $y_1,y_2\in\mR^{d_2}$,
\begin{align}\label{trans2}
\begin{split}
\p^{[\eta]}_y\bar h(y_1)&-\p^{[\eta]}_y\bar h(y_2)=\int_{\mR^{d_1}}\!\!\bigg(\big[\p_y^{[\eta]} h(x,y_1)-\p_y^{[\eta]} h(x,y_2)\big]\\
&-\sum_{\ell=1}^{[\eta]} C_{[\eta]}^\ell \cdot\frac{\p^\ell\sL_0}{\p y^\ell}(x,y_2)\big[\p_y^{[\eta]-\ell}u(x,y_1)-\p_y^{[\eta]-\ell}u(x,y_2)\big]\\
&-\sum_{\ell=1}^{[\eta]} C_{[\eta]}^\ell \Big[\frac{\p^\ell\sL_0}{\p y^\ell}(x,y_1)-\frac{\p^\ell\sL_0}{\p y^\ell}(x,y_2)\Big]\p_y^{[\eta]-\ell}u(x,y_1)\\
&-\big[\sL_0(x,y_1)-\sL_0(x,y_2)\big]\p_y^{[\eta]}u(x,y_2)\bigg)\mu^{y_1}(\dif x).
\end{split}
\end{align}
In particular, we have $\bar h\in C_b^\eta(\mR^{d_2})$.
\el

\br
(i) The above result together with Theorem \ref{popde} provide a useful tool to verify the regularity the averaged coefficients, which is a  separate problem that will always encounter in the study of averaging principle, CLT, homogenization and other limit theorems. Thus, this lemma is of independent interest.

\vspace{1mm}
(ii) Note that the left hand sides of the equality (\ref{trans}) and (\ref{trans2}) involve $\eta$-order `derivatives' of $\bar h$ with respect to the $y$-variable, while the right hand sides only involve at most $0\vee(\eta-1)$-order `derivatives' of the solution $u$ with respect to the parameter $y$. Hence these  equality can be viewed as two transfer formulas, which transfer the smoothness with respect to the parameter $y$ to the two derivatives with respect to the $x$ variable. Let us explain this more clearly when $\eta=1$. In this case, the above conclusion  simplifies to: if we have $u\in C_p^{2+\delta,0}$, then
$$
\p_y\left(\int_{\mR^{d_1}}h(x,y)\mu^y(\dif x)\right)=\int_{\mR^{d_1}}\!\!\Big[\p_y h(x,y)-\frac{\p\sL_0}{\p y}(x,y)u(x,y)\Big]\mu^y(\dif x).
$$
In order to derive the differentiablity of $\bar h(y)$, the left hand side of the above equality implies that we need to study the derivative of the invariant measure $\mu^y(\dif x)$ with respect to the parameter $y$, which is usually difficult to obtain. However, the right hand side of the above equality  only needs two derivatives of $u$ with respect to the $x$ variable (no derivative with respect to the parameter $y$ is involved), which is quite classical due to the elliptic property of the operator $\sL_0$.
\er

\begin{proof}
	We only need to prove the equality (\ref{trans}) and (\ref{trans2}). Then the assertion that $\bar h\in C_b^\eta(\mR^{d_2})$ follows directly. In fact, if (\ref{trans}) is true for $\eta\in\mN^*$, we can derive by the assumptions $a, b\in C_b^{\delta,\eta}$, $h\in C_p^{\delta,\eta}$ and  $u\in C_p^{2+\delta,(\eta-1)\vee0}$ that for some $m>0$,
	\begin{align*}
	\p_y^\eta\bar h(y)\leq C_0\int_{\mR^{d_1}}(1+|x|^m)\mu^y(\dif x)\leq C_1<\infty,
	\end{align*}
	while if (\ref{trans2}) is true for $\eta\in(0,\infty)\setminus\mN^*$, we can derive that for $y_1,y_2\in\mR^{d_2}$,
	\begin{align*}
	\big|\p^{[\eta]}_y\bar h(y_1)-\p^{[\eta]}_y\bar h(y_2)\big|&\leq C_2\Big(|y_1-y_2|^{\eta-[\eta]}\wedge1\Big)\int_{\mR^{d_1}}(1+|x|^m)\mu^{y_1}(\dif x)\\
	&\leq C_3\Big(|y_1-y_2|^{\eta-[\eta]}\wedge1\Big).
	\end{align*}
	Below, we divide the proof into two steps.
	
	\vspace{1mm}	
	\noindent(i) Let us first prove (\ref{trans}) with $\eta=1$. In fact,
	by the chain rule  we can write
	\begin{align*}
	\p_y \bar h(y)&=\p_y\left(\int_{\mR^{d_1}}h(x,y)\mu^y(\dif x)\right)\\
	&=\int_{\mR^{d_1}}\p_yh(x,y)\mu^y(\dif x)+\int_{\mR^{d_1}}h(x,y)\p_yp_\infty(x,y)\dif x.
	\end{align*}	
	Then we use a formula established in \cite[(28)]{P-V2} which yields
	\begin{align*}
	&\int_{\mR^{d_1}}\!h(x,y)\p_yp_\infty(x,y)\dif x\\
	&=-\int_{\mR^{d_1}}\!h(x,y)\left(\int_0^\infty\!\!\int_{\mR^{d_1}}p_\infty(x',y)\frac{\p\sL_0}{\p y}(x',y)p_s(x',x;y)\dif x'\dif s\right)\dif x.
	\end{align*}
	As a result, by  (\ref{pro}), Fubini's theorem and the fact that $p_t(x',x;y)$ is a density function, we  deduce that
	\begin{align*}
	&\int_{\mR^{d_1}}\!\!h(x,y)\p_yp_\infty(x,y)\dif x\\
	&=-\int_{\mR^{d_1}}\frac{\p\sL_0}{\p y}(x',y)\left(\int_0^\infty\!\!\int_{\mR^{d_1}}p_s(x',x;y)h(x,y)\dif x\dif s\right)p_\infty(x',y)\dif x'\\
	&=-\int_{\mR^{d_1}}\frac{\p\sL_0}{\p y}(x',y)\left(\int_0^\infty\!\!\int_{\mR^{d_1}}p_s(x',x;y)f(x,y)\dif x\dif s\right)p_\infty(x',y)\dif x'\\
	&=-\int_{\mR^{d_1}}\frac{\p\sL_0}{\p y}(x',y)\left(\int_0^\infty T_sf(x',y)\dif s\right)p_\infty(x',y)\dif x'\\
	&=-\int_{\mR^{d_1}}\frac{\p\sL_0}{\p y}(x',y)u(x',y)\mu^y(\dif x'),
	\end{align*}
	which in turn implies (\ref{trans}) is true for $\eta=1$.	
	The general case that $\eta\in\mN^*$ can be proved by using formula \cite[(34)]{P-V2} instead of \cite[(28)]{P-V2} and the induction argument, we omit the details here.
	
	\vspace{1mm}\noindent
	(ii) Now we prove (\ref{trans2}) when $\eta\in(0,1)$. For any $y_1, y_2\in\mR^{d_2}$, we write
	\begin{align*}
	\bar h(y_1)&-\bar h(y_2)=\int_{\mR^{d_1}}h(x,y_1)\mu^{y_1}(\dif x)-\int_{\mR^{d_1}}h(x,y_2)\mu^{y_2}(\dif x)\\
	&=\int_{\mR^{d_1}}\!\!\big[h(x,y_1)-h(x,y_2)\big]\mu^{y_1}(\dif x)+\!\int_{\mR^{d_1}}\!\!h(x,y_2)\big[\mu^{y_1}(\dif x)-\mu^{y_2}(\dif x)\big].
	\end{align*}
	Then we deduce by \cite[Lemma 4.1]{RSX} (which is similar in spirit of \cite[(28)]{P-V2}) and the same argument as above that	
	\begin{align*}
	&\int_{\mR^{d_1}}h(x,y_2)\big[\mu^{y_1}(\dif x)-\mu^{y_2}(\dif x)\big]\\
	&=-\int_0^\infty\!\!\int_{\mR^{d_1}}p_\infty(x';y_1)\big[\sL_0(x',y_1)-\sL_0(x',y_2)\big]\\
	&\qquad\qquad\qquad\times\left(\int_{\mR^{d_1}}\!p_s(x',x;y_2)h(x,y_2)\dif x\right)\!\dif x'\dif s\\
	&=-\int_0^\infty\!\!\int_{\mR^{d_1}}p_\infty(x';y_1)\big[\sL_0(x',y_1)-\sL_0(x',y_2)\big]\\
	&\qquad\qquad\qquad\times\left(\int_{\mR^{d_1}}p_s(x',x;y)f(x,y_2)\dif x\right)\dif x'\dif s\\
	&=-\int_{\mR^{d_1}}\big[\sL_0(x',y_1)-\sL_0(x',y_2)\big]u(x',y_2)\mu^{y_1}(\dif x'),
	\end{align*}
	where in the second equality we also used the fact that $p_t(x',x;y)$ is a density function. Thus, (\ref{trans2}) is true for $\eta\in(0,1)$. The general case that $\eta\in (n,n+1)$ for some $n\in\mN^*$ can be proved by using (\ref{trans}) with $\eta=n$ and the same arguments as before. The proof is finished.
\end{proof}

Now, we are in the position to give:

\begin{proof}[Proof of Theorem \ref{popde}]
	We divide the proof into four steps.
	
	\vspace{2mm}	
	\noindent{\bf Step 1.}	In this step we prove (\ref{key1}). It suffices to consider the estimate for the second order derivative $\nabla^2_xu$. To this end, we rewrite (\ref{pro}) as
	\begin{align}\label{uu}
	u(x,y)=\int_0^2T_tf(x,y)\dif t+\int_2^\infty T_tf(x,y)\dif t=:u_1(x,y)+u_2(x,y).
	\end{align}
	When $0<t\leq 2$, using (\ref{ppp1}) and the fact that $p_t(x,x';y)$ is a density function, we  deduce
	\begin{align*}
	\nabla^2_xT_tf(x,y)&=\int_{\mR^{d_1}}\nabla^2_xp_t(x,x';y)\big[f(x',y)-f(x,y)\big]\dif x'\\
	&\leq C_1[f]_{C^{\delta,0}_p}\int_{\mR^{d_1}}\big(|x-x'|^\delta\wedge1\big)\big(1+|x|^m+|x'|^m\big)\no\\
	&\qquad\qquad\qquad\qquad\times t^{-(d+2)/2}\exp\big(-c_0|x-x'|^2/t\big)\dif x'\\
	&\leq C_1[f]_{C^{\delta,0}_p}(1+|x|^m)\int_{\mR^{d_1}}\!\!t^{-(d+2-\delta)/2}\exp\big(-c_1|x-x'|^2/t\big)\dif x'\\
	&\leq C_1[f]_{C^{\delta,0}_p}(1+|x|^m)\cdot t^{(\delta-2)/2},
	\end{align*}
	while for $t>2$, we have by the semigroup property that
	\begin{align*}
	&\nabla^2_xT_tf(x,y)=\int_{\mR^{d_1}}\int_{\mR^{d_1}}\nabla^2_xp_1(x,z;y)p_{t-1}(z,x';y)\dif zf(x',y)\dif x'\\
	&= \int_{\mR^{d_1}}\int_{\mR^{d_1}}\nabla^2_xp_1(x,z;y)\big[p_{t-1}(z,x';y)-p_\infty(x';y)\big]\dif zf(x',y)\dif x'\\
	&\leq C_2[f]_{C^{\delta,0}_p}\int_{\mR^{d_1}}\int_{\mR^{d_1}}\exp\big(-c_0|x-z|^2\big)\frac{(1+|z|^m)(1+|x'|^m)}{(1+t)^k(1+|x'|^j)}\dif z\dif x'\no\\
	&\leq C_2[f]_{C^{\delta,0}_p}\frac{(1+|x|^m)}{(1+t)^k}\int_{\mR^{d_1}}\frac{(1+|x'|^m)}{(1+|x'|^j)}\dif x'\\
	&\leq C_2[f]_{C^{\delta,0}_p}\frac{(1+|x|^m)}{(1+t)^k},
	\end{align*}
	where in the third inequality we have used (\ref{pin2}), and we choose $j>d_1+m$ in the last inequality. Now, taking these two estimates back into (\ref{uu}) gives
	\begin{align*}
	|\nabla^2_xu(x,y)|&\leq C_3[f]_{C^{\delta,0}_p}(1+|x|^m)\bigg(\int_0^2t^{(\delta-2)/2}\dif t+\int_2^\infty \frac{1}{1+t^k}\dif t\bigg)\\
	&\leq C_3[f]_{C^{\delta,0}_p}(1+|x|^m).
	\end{align*}
	Thus estimate (\ref{key1}) is true.
	
	\vspace{2mm}	
	\noindent{\bf Step 2.} We proceed to prove estimate (\ref{key11}). Note that when $|x_1-x_2|\geq 1$, the conclusion follows directly from (\ref{key1}). Below, we focus on the case where $|x_1-x_2|<1$. Let $\tilde x$ be one of the two points $x_1$ and $x_2$ which is nearer to $x'$, and denote by $\cO$ the ball with center $\tilde x$ and radius $2|x_1-x_2|$. Without loss of generality, we may assume that $\tilde x=x_1$. We write
	\begin{align*}
	\nabla&^2_xu_1(x_1,y)-\nabla^2_xu_1(x_2,y)=\!\int_0^2\!\!\int_{\cO}\nabla^2_x p_t(x_1,x';y)\big[f(x',y)-f(x_1,y)\big]\dif x'\dif t\\
	&-\int_0^2\!\!\int_{\cO}\nabla^2_xp_t(x_2,x';y)\big[f(x',y)-f(x_2,y)\big]\dif x'\dif t\\
	&+\int_0^2\!\!\int_{\mR^{d_1}\setminus\cO}\!\big[\nabla^2_x p_t(x_1,x';y)-\!\nabla^2_x p_t(x_2,x';y)\big]\big[f(x',y)-f(x_1,y)\big]\dif x'\dif t\\
	&+\big[f(x_2,y)-f(x_1,y)\big]\int_0^2\!\!\int_{\mR^{d_1}\setminus\cO}\nabla^2_x p_t(x_2,x';y)\dif x'\dif t=:\sum_{i=1}^4\cI_i.
	\end{align*}
	For the first term, we have by (\ref{ppp1}) that
	\begin{align*}
	\cI_1&\leq C_1[f]_{C^{\delta,0}_p}\!\!\int_0^2\!\!\!\int_{\cO}t^{-(d+2)/2}\exp\!\big(\!-c_0|x_1-x'|^2/t\big)\\
	&\qquad\qquad\qquad\qquad\times|x'-x_1|^\delta\big(1+|x_1|^m+|x'|^m\big)\dif x'\dif t\\
	&\leq C_1[f]_{C^{\delta,0}_p}(1+|x_1|^m)\int_{\cO}\int_0^2t^{-(d+2-\delta)/2}\exp\big(-c_1|x_1-x'|^2/t\big)\dif t\dif x'\\
	&\leq C_1[f]_{C^{\delta,0}_p}(1+|x_1|^m)\int_{\cO}|x_1-x'|^{-d+\delta}\dif x'\\
	&\leq C_1[f]_{C^{\delta,0}_p}|x_1-x_2|^\delta(1+|x_1|^m).
	\end{align*}
	In completely the same way we get
	\begin{align*}
	\cI_2&\leq C_2[f]_{C^{\delta,0}_p}|x_1-x_2|^\delta(1+|x_2|^m).
	\end{align*}
	To control the third term, note that for $x'\in \mR^{d_1}\setminus\cO$, we have
	$$
	|x_2-x'|/2\leq|x_1-x'|\leq 3|x_2-x'|/2.
	$$
	As a result, by (\ref{ppp2}) we deduce  that
	\begin{align*}
	\cI_3&\leq C_3[f]_{C^{\delta,0}_p}|x_1-x_2|\int_0^2\!\!\int_{\mR^{d_1}\setminus\cO}t^{-(d+3)/2}\exp\big(-c_0|x_1-x'|^2/t\big)\\
	&\qquad\qquad\qquad\qquad\qquad\times|x_1-x'|^\delta\big(1+|x_1|^m+|x'|^m\big)\dif x'\dif t\\
	&\quad+C_3[f]_{C^{\delta,0}_p}|x_1-x_2|^\beta\int_0^2\!\!\int_{\mR^{d_1}\setminus\cO}t^{-(d+2)/2}\exp\big(-c_0|x_1-x'|^2/t\big)\\
	&\qquad\qquad\qquad\qquad\times|x_1-x'|^\delta\big(1+|x_1|^m+|x'|^m\big)\dif x'\dif t=:\cI_{31}+\cI_{32}.
	\end{align*}
	We further control $\cI_{31}$ by
	\begin{align*}
	\cI_{31}&\leq C_3[f]_{C^{\delta,0}_p}|x_1-x_2|(1+|x_1|^m)\\
	&\quad\times\int_{\mR^{d_1}\setminus\cO}\!\!\left(\int_0^2t^{-(d+3)/2}\exp\big(-c_1|x_1-x'|^2/t\big)\dif t\right)|x_1-x'|^\delta\dif x'\\
	&\leq C_3[f]_{C^{\delta,0}_p}|x_1-x_2|(1+|x_1|^m)\int_{\mR^{d_1}\setminus\cO}|x-x'|^{-d-1+\delta}\dif x'\\
	&\leq C_3[f]_{C^{\delta,0}_p}|x_1-x_2|^\delta(1+|x_1|^m),
	\end{align*}
	and it is easy to check that
	\begin{align*}
	\cI_{32}&\leq C_3[f]_{C^{\delta,0}_p}|x_1-x_2|^\delta(1+|x_1|^m)\\
	&\quad\times\int_0^2\!\!\int_{\mR^{d_1}}t^{-(d+2)/2}\exp\big(-c_2|x_1-x'|^2/t\big)|x_1-x'|^\delta\dif x'\dif t\\
	&\leq C_3[f]_{C^{\delta,0}_p}|x_1-x_2|^\delta(1+|x_1|^m)\int_0^2t^{\delta/2-1}\dif t\\
	&\leq C_3[f]_{C^{\delta,0}_p}|x_1-x_2|^\delta(1+|x_1|^m).
	\end{align*}
	Finally, thanks to (\ref{ppp3}), we have
	\begin{align*}
	\cI_4&\leq C_4[f]_{C^{\delta,0}_p}|x_1-x_2|^\delta(1+|x_1|^m+|x_2|^m)\int_0^2t^{\delta/2-1}\dif t\\
	&\leq C_4[f]_{C^{\delta,0}_p}|x_1-x_2|^\delta(1+|x_1|^m+|x_2|^m).
	\end{align*}
	Combining the above computations, we arrive at
	$$
	|\nabla^2_xu_1(x_1,y)-\nabla^2_xu_1(x_2,y)|\leq C_5[f]_{C^{\delta,0}_p}|x_1-x_2|^\delta(1+|x_1|^m+|x_2|^m).
	$$
	On the other hand, we  derive as in Step 1 that for $t\geq 2$,
	\begin{align*}
	&\nabla^2_xT_tf(x_1,y)-\nabla^2_xT_tf(x_2,y)\\
	&=\int_{\mR^{d_1}}\int_{\mR^{d_1}}\big[\nabla^2_xp_1(x_1,z;y)-\nabla^2_xp_1(x_2,z;y)\big]p_{t-1}(z,x';y)\dif zf(x',y)\dif x'\\
	&= \int_{\mR^{d_1}}\int_{\mR^{d_1}}\!\big[\nabla^2_xp_1(x_1,z;y)\!-\!\nabla^2_xp_1(x_2,z;y)\big]\\
	&\qquad\qquad\qquad\qquad\qquad\times\big[p_{t-1}(z,x';y)\!-\!p_\infty(x';y)\big]\dif zf(x',y)\dif x'\\
	&\leq C_6[f]_{C^{\delta,0}_p}|x_1-x_2|^\delta\int_{\mR^{d_1}}\int_{\mR^{d_1}}\exp\big(-c_0|x_1-z|^2\big)\\
	&\qquad\qquad\qquad\qquad\qquad\qquad\times\frac{(1+|z|^m)(1+|x'|^m)}{(1+t)^k(1+|x'|^j)}\dif z\dif x'\\
	&\leq C_6[f]_{C^{\delta,0}_p}|x_1-x_2|^\delta\frac{(1+|x_1|^m)}{(1+t)^k},
	\end{align*}
	which in turn yields
	\begin{align*}
	|\nabla^2_xu_2(x_1,y)-\nabla^2_xu_2(x_2,y)|&\leq C_6[f]_{C^{\delta,0}_p}|x_1-x_2|^\delta\int_2^\infty\frac{(1+|x_1|^m)}{(1+t)^k}\dif t\\
	&\leq C_6[f]_{C^{\delta,0}_p}|x_1-x_2|^\delta(1+|x_1|^m).
	\end{align*}
	Now using (\ref{uu}) again, we   get (\ref{key11}).
	
	\vspace{2mm}	
	\noindent{\bf Step 3.} In this step, we prove estimate (\ref{key3}) when $\eta\in\mN^*$. We shall only focus on the a priori estimates. Let us first consider the case $\eta=1$. We start form the equation itself, i.e., $u$ is a classical solution to
	$$
	\sL_0(x,y)u(x,y)=f(x,y).
	$$
	Taking partial derivative with respect to the $y$ variable from both sides of the equation, we  get that
	\begin{align}\label{ny1}
	\sL_0(x,y)\p_y u(x,y)=\p_y f(x,y)-\frac{\p\sL_0}{\p y}(x,y)u(x,y).
	\end{align}
	According to (\ref{trans}), the right hand side of (\ref{ny1}) satisfies the centering condition (\ref{cen2}), i.e.,
	\begin{align*}
	\int_{\mR^{d_1}}\!\!\Big[\p_y  f(x,y)-\frac{\p\sL_0}{\p y}(x,y)u(x,y)\Big]\mu^y(\dif x)=\p_y\left(\int_{\mR^{d_1}}f(x,y)\mu^y(\dif x)\right)=0.
	\end{align*}
	Moreover, by assumption we have that $\p_y f\in C^{\delta,0}_p$. Meanwhile, due to the assumptions $a,b\in C_b^{\delta,1}$ and in view of the a priori estimates (\ref{key1})-(\ref{key11}), it is easily checked that for any $x_1,x_2\in\mR^{d_1}$ and $y\in\mR^{d_2}$, there exists a constant $C_1>0$ such that
	\begin{align*}
	&\Big|\frac{\p\sL_0}{\p y}(x_1,y)u(x_1,y)-\frac{\p\sL_0}{\p y}(x_2,y)u(x_2,y)\Big|\\
	&\leq C_1\Big(\|a\|_{C_b^{\delta,1}}+\|b\|_{C_b^{\delta,1}}\Big)\|f\|_{C^{\delta,0}_p}\big(|x_1-x_2|^\delta\wedge1\big)(1+|x_1|^m+|x_2|^m).
	\end{align*}
	As a result, we have $\tfrac{\p\sL_0}{\p y}u\in C_p^{\delta,0}$ with
	$$
	\big[\tfrac{\p\sL_0}{\p y}u\big]_{C_p^{\delta,0}}\leq C_1\Big(\|a\|_{C_b^{\delta,1}}+\|b\|_{C_b^{\delta,1}}\Big)\|f\|_{C^{\delta,0}_p}.
	$$
	Thus, using the conclusions proved in Step 1 for the function $\p_yu$, we   get that $\p_yu\in C_p^{2+\delta,0}$ and by the a priori estimate (\ref{key1}), we have
	\begin{align*}
	|\p_yu(x,y)|&\leq C_2\Big([\p_yf]_{C^{\delta,0}_p}+\big[\tfrac{\p\sL_0}{\p y}u\big]_{C_p^{\delta,0}}\Big)(1+|x|^m),\\
	&\leq C_2\Big([f]_{C^{\delta,1}_p}+\|a\|_{C_b^{\delta,1}}+\|b\|_{C_b^{\delta,1}}\Big)(1+|x|^m),
	\end{align*}
	which in particular yields (\ref{key3}) for $\eta=1$. Suppose  that (\ref{key3}) holds for some $\eta=n$ with constant $\cK_\eta$ given by (\ref{dep}). By induction we  find that for $\eta=n+1$,
	\begin{align*}
	\sL_0(x,y)\p_y^{n+1} u(x,y)=\p_y^{n+1} f(x,y)-\sum_{\ell=1}^{n+1} C_{n+1}^\ell \cdot\frac{\p^\ell\sL_0}{\p y^\ell}(x,y)\p_y^{n+1-\ell}u(x,y).
	\end{align*}
	According to (\ref{trans}), the right hand side of the above equality satisfies the centering condition (\ref{cen2}), i.e.,
	\begin{align*}
	&\int_{\mR^{d_1}}\!\!\Big[\p_y^{n+1} f(x,y)-\sum_{\ell=1}^{n+1} C_{n+1}^\ell \cdot\frac{\p^\ell\sL_0}{\p y^\ell}(x,y)\p_y^{n+1-\ell}u(x,y)\Big]\mu^y(\dif x)\\
	&=\p^{n+1}_y\left(\int_{\mR^{d_1}}f(x,y)\mu^y(\dif x)\right)=0.
	\end{align*}
	It then follows by using (\ref{key1}) again for the function $\p^{n+1}_yu$ and the induction assumption that
	\begin{align*}
	|\p_y^{n+1} u(x,y)|
	&\leq C_3\Big([\p^{n+1}_yf]_{C_p^{\delta,0}}+\sum_{\ell=1}^{n+1} C_{n+1}^\ell \cdot\kappa_\ell\big[\p^{n+1-\ell}_yu\big]_{C_p^{\delta,0}}\Big)(1+|x|^m)\\
	&\leq C_3\Big([f]_{C_p^{\delta,n+1}}+\sum_{\ell=1}^{n+1} C_{n+1}^\ell \cdot\kappa_\ell\cdot\cK_{n+1-\ell}\Big)(1+|x|^m),
	\end{align*}
	which means (\ref{key3}) is true for $\eta=n+1$.
	
	\vspace{2mm}	
	\noindent{\bf Step 4.}  Finally, we prove (\ref{key3}) when $\eta\in (0,\infty)\setminus\mN^*$. Let us first consider the case $\eta\in(0,1)$. By the equation (\ref{pde1}), we write for any $x\in\mR^{d_1}$ and $y_1,y_2\in\mR^{d_2}$
	\begin{align}\label{Holl1}
	\sL_0(x,y_1)\big[u(x,y_1)-u(x,y_2)\big]&=\big[f(x,y_1)-f(x,y_2)\big]\no\\
	&\quad+\big[\sL_0(x,y_2)-\sL_0(x,y_1)\big]u(x,y_2).
	\end{align}
	Using (\ref{trans2}) with $\eta\in(0,1)$ (thus $[\eta]=0$ and the sum terms in (\ref{trans2}) do not appear), we obtain for the right hand side of (\ref{Holl1}) that
	\begin{align*}
	&\int_{\mR^{d_1}}\Big(\big[f(x,y_1)-f(x,y_2)\big]+\big[\sL_0(x,y_2)-\sL_0(x,y_1)\big]u(x,y_2)\Big)\mu^{y_1}(\dif x)\\
	&=\int_{\mR^{d_1}}f(x,y_1)\mu^{y_1}(\dif x)-\int_{\mR^{d_1}}f(x,y_2)\mu^{y_2}(\dif x)=0.
	\end{align*}
	As a result, we can use (\ref{key1}) for the function $u(x,y_1)-u(x,y_2)$ and by the same argument as above to get that
	\begin{align*}
	\big|u(x,y_1)-u(x,y_2)\big|&\leq C_1\big(|y_1-y_2|^\eta\wedge1\big)\\
	&\quad\quad\times\Big([f]_{C_p^{\delta,\eta}}+\|a\|_{C_b^{\delta,\eta}}+\|b\|_{C_b^{\delta,\eta}}\Big)(1+|x|^m),
	\end{align*}
	which means that $u(x,\cdot)\in C_b^\eta$ and (\ref{key3}) is true. Assume that (\ref{key3}) holds for some $\eta\in(n,n+1)$ with constant $\cK_\eta$ given by (\ref{dep2}). Then  for $\eta\in(n+1,n+2)$, by the conclusion proved in (ii),  for any $x\in\mR^{d_1}$ and $y_1,y_2\in\mR^{d_2}$ we obtain
	\begin{align}\label{Holn}
	\sL_0(x,y_1)&\big[\p_y^{n+1} u(x,y_1)-\p_y^{n+1} u(x,y_2)\big]=\big[\p_y^{n+1} f(x,y_1)-\p_y^{n+1} f(x,y_2)\big]\no\\
	&-\sum_{\ell=1}^{n+1} C_{n+1}^\ell \cdot\frac{\p^\ell\sL_0}{\p y^\ell}(x,y_2)\big[\p_y^{n+1-\ell}u(x,y_1)-\p_y^{n+1-\ell}u(x,y_2)\big]\no\\
	&-\sum_{\ell=1}^{n+1} C_{n+1}^\ell \Big[\frac{\p^\ell\sL_0}{\p y^\ell}(x,y_1)-\frac{\p^\ell\sL_0}{\p y^\ell}(x,y_2)\Big]\p_y^{n+1-\ell}u(x,y_1)\no\\
	&-\big[\sL_0(x,y_1)-\sL_0(x,y_2)\big]\p_y^{n+1} u(x,y_2).
	\end{align}
	Using (\ref{trans2}) again  we can see that the right hand side of (\ref{Holn}) satisfies the centering condition. Consequently, we have
	\begin{align*}
	&\big|\p_y^{n+1} u(x,y_1)-\p_y^{n+1} u(x,y_2)\big|\\
	&\leq C_2\big(|y_1-y_2|^\eta\wedge1\big)\Big([f]_{C_p^{\delta,\eta}}+\kappa_{\eta-n-1}[\p^{n+1}_yu]_{C_p^{\delta,0}}\\
	&+\sum_{\ell=1}^{n+1} C_{n+1}^\ell\cdot\kappa_\ell[\p^{\eta-\ell}_yu]_{C_p^{\delta,0}}+\sum_{\ell=1}^{n+1} C_{n+1}^\ell\cdot\kappa_{\eta+\ell-n-1}[\p^{n+1-\ell}_yu]_{C_p^{\delta,0}}\Big)(1+|x|^m)\\
	&\leq C_2\big(|y_1-y_2|^\eta\wedge1\big)\Big([f]_{C_p^{\delta,\eta}}+\kappa_{\eta-n-1}\cK_{n+1}\\
	&+\sum_{\ell=1}^{n+1} C_{n+1}^\ell\cdot\kappa_\ell\cK_{\eta-\ell}+\sum_{\ell=1}^{n+1} C_{n+1}^\ell\cdot\kappa_{\eta+\ell-n-1}\cK_{n+1-\ell}\Big)(1+|x|^m)\\
	&=C_2\cK_\eta\big(|y_1-y_2|^\eta\wedge1\big)(1+|x|^m),
	\end{align*}
	which means (\ref{key3}) is true for $\eta\in (n+1,n+2)$. So, the proof is finished.
\end{proof}

\section{Fluctuation estimates}

We shall first prepare some results concerning the mollifying approximation of functions in Subsection 4.1. Then by using the Poisson equation (\ref{pde0}) we derive two new fluctuation estimates for SDE (\ref{sde0}): one of functional LLN type  in Subsection 4.2 and one of functional CLT type in Subsection 4.3.

\subsection{Mollifying approximation}
We  need some mollification arguments due to our low regularity assumptions on the coefficients. To this end, let $\rho_1:\mR\to[0,1]$ and $\rho_2:\mR^{d_2}\to[0,1]$ be two smooth radial convolution kernel functions
such that
$\int_\mR\rho_1(r)\dif r=\int_{\mR^{d_2}}\rho_2(y)\dif y=1$, and for any $k\geq 1$, there exist constants $C_k>0$ such that $|\nabla^k\rho_1(r)|\leq C_k\rho_1(r)$ and  $|\nabla^k\rho_2(y)|\leq C_k\rho_2(y)$. For every $n\in\mN^*$, set
$$
\rho_1^n(y):=n^2\rho_1(n^2r)\quad \text{and}\quad
\rho_2^n(y):=n^{d_2}\rho_2(ny).$$
Note that the scaling speed of $\rho_2^n$ is $n$, while $\rho_1^n$ has the speed $n^2$ in mollifying.
Given a function $f(t,x,y)$, define the mollifying approximations of $f$ in $t$ and $y$ variables by
\begin{align}\label{fn}
f_n(t,x,y):=f*\rho_2^n*\rho_1^n:=\int_{\mR^{d_2+1}}f(t-s,x,y-z)\rho_2^{n}(z)\rho_1^n(s)\dif z\dif s.
\end{align}
We have the following easy result.

\bl
Let $f\in C_p^{\vartheta/2,\delta,\vartheta}$ with $0<\vartheta\leq 2$, $0<\delta\leq 1$, and define $f_n$ by (\ref{fn}). Then we have
\begin{align}\label{n111}
\|f(\cdot,x,\cdot)-f_n(\cdot,x,\cdot)\|_\infty\leq C_0n^{-\vartheta}(1+|x|^m),
\end{align}
and
\begin{align}\label{n222}
\|\p_tf_n(\cdot,x,\cdot)\|_{\infty}+\|\nabla^2_yf_n(\cdot,x,\cdot)\|_{\infty}\leq C_0n^{2-\vartheta}(1+|x|^m),
\end{align}
where $C_0>0$ is a constant independent of $n$.
\el
\begin{proof}
	We divide the proof into two cases.\\
	(i)	(Case $0<\vartheta\leq 1$).
	By definition and a change of variable, there exists a constant $m>0$ such that
	\begin{align*}
	|f(t,x,y)\!-\!f_n(t,x,y)|&\leq \!\int_{\mR^{d_2+1}}\!\!\big|f(t,x,y)-f(t-s,x,y-z)\big|\rho_2^{n}(z)\rho_1^n(s)\dif z\dif s\\
	&\leq C_0\int_{\mR^{d_2+1}}\big(s^{\vartheta/2}+|z|^{\vartheta}\big)(1+|x|^m)\rho_2^{n}(z)\rho_1^n(s)\dif z\dif s\\
	&\leq C_0n^{-\vartheta}(1+|x|^m).
	\end{align*}
	Furthermore,
	\begin{align*}
	&|\p_tf_n(t,x,y)|\leq \int_{\mR^{d_2+1}}\!\!\big|f(t-s,x,y-z)-f(t,x,y-z)\big|\rho_2^{n}(z)|\p_s\rho_1^n(s)|\dif z\dif s\\
	&\leq C_0n^2\int_{\mR^{d_2+1}}\!\!s^{\vartheta/2}(1+|x|^m)\rho_2^{n}(z)\rho_1^n(s)\dif z\dif s\leq C_0n^{2-\vartheta}(1+|x|^m),
	\end{align*}
	and
	\begin{align*}
	&|\nabla_y^2f_n(t,x,y)|\!\leq\!\! \int_{\mR^{d_2+1}}\!\!\big|f(t-s,x,y-z)-f(t-s,x,y)\big| |\nabla_z^2\rho_2^{n}(z)|\rho_1^n(s)\dif z\dif s\\
	&\leq C_0n^2\int_{\mR^{d_2+1}}|z|^{\vartheta}(1+|x|^m)\rho_2^{n}(z)\rho_1^n(s)\dif z\dif s\leq C_0n^{2-\vartheta}(1+|x|^m).
	\end{align*}
	(ii)	(Case $1<\vartheta\leq 2$). By the symmetric property of $\rho_1$ and $\rho_2$, we can  use the second order difference to deduce that
	\begin{align*}
	|f(t,x,y)-f_n(t,x,y)|&\leq \int_{\mR^{d_2+1}}\big|f(t-s,x,y-z)+f(t-s,x,y+z)\\
	&\quad\quad\qquad\quad-2f(t,x,y)\big|\cdot\rho_2^{n}(z)\rho_1^n(s)\dif z\dif s\\
	&\leq C_1\int_{\mR^{d_2+1}}\big(s^{\vartheta/2}+|z|^{\vartheta}\big)(1+|x|^m)\cdot\rho_2^{n}(z)\rho_1^n(s)\dif z\dif s\\
	&\leq C_1n^{-\vartheta}(1+|x|^m),
	\end{align*}
	and
	\begin{align*}
	|\nabla_y^2f_n(t,x,y)|&\leq \int_{\mR^{d_2+1}}\!\big|\nabla_yf(t-s,x,y-z)\\
	&\quad\quad\qquad\quad-\nabla_yf(t-s,x,y)\big||\nabla_z\rho_2^{n}(z)|\rho_1^n(s)\dif z\dif s\\
	&\leq C_2n\int_{\mR^{d_2+1}}|z|^{\vartheta-1}(1+|x|^m)\cdot\rho_2^{n}(z)\rho_1^n(s)\dif z\dif s\\
	&\leq C_2n^{2-\vartheta}(1+|x|^m).
	\end{align*}
	The estimate concerning the derivative with respect to the $t$ variable can be proved similarly. The proof is finished.
\end{proof}

\subsection{Fluctuation estimate - LLN type}

Given a function $h(t,x,y)$, recall that  $f(t,x,y):=h(t,x,y)-\bar h(t,y)$ satisfies the centering condition (\ref{cen2}), where $\bar h$ is defined by (\ref{aver}).
The following result establishes the behavior of the  fluctuation between $h(s,X_s^\eps,Y_s^\eps)$ and $\bar h(s, Y_s^\eps)$ over the time interval $[0,t]$, which will play an important role below.

\bl\label{key}
Let {\bf (A$_\sigma$)}, {\bf (A$_b$)} and (\ref{ac}) hold  true. Assume that $b,\sigma\in C_b^{\delta,\vartheta}$  with $0<\delta,\vartheta\leq 2$, and $c, F, H, G\in L^\infty_p$. Then for every $f\in C_p^{\vartheta/2,\delta,\vartheta}$  satisfying (\ref{cen2}), we have
\begin{align*}
\mE\left(\int_0^tf(s,X_s^\eps,Y_s^\eps)\dif s\right)\leq C_t\Big(\alpha_\eps^{\vartheta}+\alpha_\eps^{\vartheta\wedge1}\cdot\frac{\alpha_\eps}{\gamma_\eps}+\frac{\alpha_\eps^2}{\beta_\eps}\Big),
\end{align*}
where $C_t>0$ is a constant independent of $\delta,\eps$.
\el

\br
If  $0<\vartheta\leq1$, then the fluctuation is controlled by $\alpha_\eps^{\vartheta}+\frac{\alpha_\eps^2}{\beta_\eps}$; and if $\vartheta=2$, the error bound can be controlled by $\frac{\alpha_\eps^2}{\gamma_\eps}+\frac{\alpha_\eps^2}{\beta_\eps}$; while when $1<\vartheta\leq 2$, the fluctuation will depend on the balance between $\alpha_\eps$ and $\gamma_\eps$.
In particular, in the case of Regime 3 and Regime 4, the bound simplifies to $\alpha_\eps^{\vartheta\wedge1}$.
\er
\begin{proof}
	Since $f$ satisfies (\ref{cen2}), by Theorem \ref{popde} and in view of (\ref{pro}), there exists a unique solution $\Phi^f(t,x,y)\in C^{\vartheta/2,2+\delta,\vartheta}_p$ to the following Poisson equation in $\mR^{d_1}$:
	\begin{align}\label{pn}
	\sL_0(x,y)\Phi^f(t,x,y)=-f(t,x,y),
	\end{align}
	where $(t,y)\in\mR_+\times\mR^{d_2}$ are parameters.
	Let $\Phi_n^f$ be the mollifyer  of $\Phi^f$ defined as in (\ref{fn}).
	By It\^o's formula, we  get
	\begin{align*} \Phi_n^f(t,X_{t}^\eps,Y_{t}^\eps)&=\Phi_n^f(0,x,y)+\frac{1}{\alpha_\eps}M^1_{n}(t)+M^2_{n}(t)\\
	&\quad+\int_0^{t}\Big(\p_s+\beta_\eps^{-1}\sL_3+\gamma_\eps^{-1}\sL_2+\sL_1\Big)\Phi_n^f(s,X_s^\eps,Y_s^\eps)\dif s\\
	&\quad+\frac{1}{{\alpha_\eps^2}}\int_0^{t}\sL_0\Phi_n^f(s,X_s^\eps,Y_s^\eps)\dif s,
	\end{align*}
	where $\sL_3, \sL_2$ and $\sL_1$ are defined by (\ref{lll}), and for $i=1,2$, $M^i_{n}(t)$ are martingales given by
	$$
	M^1_{n}(t):=\int_0^t\nabla_x\Phi_n^f(s,X_s^\eps,Y_s^\eps)\sigma(X_s^\eps,Y_s^\eps)\dif W_s^1
	$$
	and
	$$
	M^2_{n}(t):=\int_0^t\nabla_y\Phi_n^f(s,X_s^\eps,Y_s^\eps)G(s,X_s^\eps,Y_s^\eps)\dif W_s^2.
	$$
	This in turn yields that
	\begin{align}\label{ito}
	\int^{t}_0\!f(s,X^{\eps}_s, Y^{\eps}_s)\dif s
	&={\alpha_\eps^2}\Phi_n^f(0,x,y)-{\alpha_\eps^2}\Phi_n^f({t},X_{t}^\eps,Y_{t}^\eps)+{\alpha_\eps} M^1_{n}(t)\no\\
	&+{\alpha_\eps^2}M^2_{n}(t)+{\alpha_\eps^2}\int_0^{t}\big(\p_s+\sL_1\big)\Phi_n^f(s,X_s^\eps,Y_s^\eps)\dif s\no\\
	&+\frac{\alpha_\eps^2}{\gamma_\eps}\int_0^{t}\!\sL_2\Phi_n^f(s,X_s^\eps,Y_s^\eps)\dif s+\frac{\alpha_\eps^2}{\beta_\eps}\int_0^{t}\!\sL_3\Phi_n^f(s,X_s^\eps,Y_s^\eps)\dif s\no\\
	&+\int_0^{t}\big(\sL_0\Phi^f_n-\sL_0\Phi^f\big)(s,X_s^\eps,Y_s^\eps)\dif s.
	\end{align}
	As a result, we have
	\begin{align*}
	\cQ(\eps)&:=\mE\left(\int_0^tf(s,X_s^\eps,Y_s^\eps)\dif s\right)\leq2{\alpha_\eps^2}\mE\|\Phi^f_n(\cdot,X_t^\eps,\cdot)\|_\infty\\
	&\quad+{\alpha_\eps^2}\mE\left|\int_0^{t}\big(\p_s+\sL_1\big)\Phi_n^f(s,X_s^\eps,Y_s^\eps)\dif s\right|\\
	&\quad+\frac{\alpha_\eps^2}{\gamma_\eps}\mE\left|\int_0^{t}\sL_2\Phi_n^f(s,X_s^\eps,Y_s^\eps)\dif s\right|+\frac{\alpha_\eps^2}{\beta_\eps}\mE\left|\int_0^{t}\sL_3\Phi_n^f(s,X_s^\eps,Y_s^\eps)\dif s\right|\\
	&\quad+\mE\left|\int_0^{t}\big(\sL_0\Phi_n^f-\sL_0\Phi^f\big)(s,X_s^\eps,Y_s^\eps)\dif s\right|=:\sum_{i=1}^5\cQ_{i}(\eps).
	\end{align*}
	Note that the assumptions {\bf (A$_\sigma$)} and (\ref{ac}) hold uniformly in $y$. Hence it follows by \cite[Lemma 1]{Ve5} (see also \cite[Lemma 2]{P-V2} or \cite{S1}) that for any $m>0$ and $\eps\in(0,1/2)$,
	\begin{align}
	\mE|X^{\eps}_t|^{m}\leq C_0(1+|x|^{m}),\label{ME}
	\end{align}
	where $C_0$ is a positive constant independent of $\eps$.
	Hence, it follows by (\ref{key1}) that there exists a constant $C_1>0$ such that
	$$
	\cQ_{1}(\eps)\leq C_1\alpha_\eps^2\mE\|\Phi^f(\cdot,X_t^\eps,\cdot)\|_\infty\leq C_1\alpha_\eps^2\mE\big(1+|X_t^\eps|^m\big)\leq C_1\alpha_\eps^2.
	$$	
	To control the second term, we have by  (\ref{n222}) and the assumption $F, G\in L^\infty_p$ that
	\begin{align*}
	&\|\big(\p_s+\sL_1\big)\Phi_n^f(\cdot,x,\cdot)\|_\infty\leq C_2\big(1+|x|^m\big)\\
	&\qquad\times\Big(\|\p_s\Phi_n^f(\cdot,x,\cdot)\|_\infty +\sum_{\ell=1,2}\big\|\nabla_y^\ell\Phi_n^f(\cdot,x,\cdot)\big\|_\infty\Big)\leq C_2n^{2-\vartheta}(1+|x|^{2m}).
	\end{align*}
	Consequently, by (\ref{ME})
	$$
	\cQ_{2}(\eps)\leq C_2{\alpha_\eps^2}n^{2-\vartheta}\mE\left(\int_0^t\big(1+|X_s^\eps|^{2m}\big)\dif s\right)\leq C_2{\alpha_\eps^2}n^{2-\vartheta}.
	$$
	Following the same argument as above, we  get that when $0<\vartheta\leq 1$,
	\begin{align*}
	\|\sL_2\Phi_n^f(\cdot,x,\cdot)\|_\infty\leq C_3\|\nabla_y\Phi^f_n(\cdot,x,\cdot)\|_\infty(1+|x|^m)\leq C_3 n^{1-\vartheta}(1+|x|^{2m}),
	\end{align*}
	while for $1<\vartheta\leq 2$, we have
	\begin{align*}
	\|\sL_2\Phi_n^f(\cdot,x,\cdot)\|_\infty\leq C_3\|\nabla_y\Phi^f(\cdot,x,\cdot)\|_\infty(1+|x|^m)\leq  C_3 (1+|x|^{2m}).
	\end{align*}
	Thus, we  get
	\begin{align*}
	\cQ_{3}(\eps)
	&\leq C_3\frac{\alpha_\eps^2}{\gamma_\eps} n^{1-(\vartheta\wedge1)}\mE\left(\int_0^t\big(1+|X_s^\eps|^m\big)\dif s\right)\leq C_3\frac{\alpha_\eps^2}{\gamma_\eps} n^{1-(\vartheta\wedge1)}.
	\end{align*}
	Meanwhile, since $c\in L^\infty_p$, it follows easily that for some $m>0$,
	$$
	\cQ_{4}(\eps)\leq C_4\frac{\alpha_\eps^2}{\beta_\eps}\mE\left(\int_0^t\big(1+|X_s^\eps|^m\big)\dif s\right)\leq C_4\frac{\alpha_\eps^2}{\beta_\eps}.
	$$
	Finally, since $\nabla_x^2\Phi^f\in C_p^{\vartheta/2,\delta,\vartheta}$ and by the fact that
	$$
	\nabla_x^2(\Phi^f_n)=(\nabla_x^2\Phi^f)*\rho_1^n*\rho_2^n,
	$$
	we can derive by (\ref{n111}) that
	\begin{align*}
	\cQ_{5}(\eps)
	&\leq C_5\big(\|a\|_\infty+\|b\|_\infty\big)\cdot\mE\left(\int_0^t\!\sum_{\ell=1,2}\big\|(\nabla_x^\ell\Phi_n^f-\nabla_x^\ell\Phi^f)(\cdot,X_s^\eps,\cdot)\big\|_\infty\dif s\right)\\
	&\leq C_5n^{-\vartheta}\mE\left(\int_0^t\big(1+|X_s^\eps|^m\big)\dif s\right)\leq C_5 n^{-\vartheta}.
	\end{align*}
	Combining the above computations, we arrive at
	$$
	\cQ(\eps)\leq C_6\Big({\alpha_\eps^2}+{\alpha_\eps^2}n^{2-\vartheta}+\frac{\alpha_\eps^2}{\gamma_\eps} n^{1-(\vartheta\wedge1)}+n^{-\vartheta}+\frac{\alpha_\eps^2}{\beta_\eps}\Big).
	$$	
	Taking $n={\alpha_\eps^{-1}}$, we thus get
	$$
	\cQ(\eps)\leq C_6\Big(\alpha_\eps^\vartheta+\frac{\alpha_\eps^{1+(\vartheta\wedge1)}}{\gamma_\eps} +\frac{\alpha_\eps^2}{\beta_\eps}\Big).
	$$
	The proof is finished.
\end{proof}

\subsection{Fluctuation estimate - CLT type}

Now, we derive a CLT type fluctuation estimate for $f(s,X_s^\eps,Y_s^\eps)$ over the time interval $[0,t]$. This will depend on the orders how $\alpha_\eps, \beta_\eps, \gamma_\eps$ go to zero described in Regime 1-Regime 4 of (\ref{regime}), and  the limit behavior will involve the solution of the auxiliary Poisson equation.

Recall that $\Phi^f$ is the solution to the Poisson equation (\ref{pn}).
For simplify, we denote by
\begin{align*}
\overline{c\cdot\nabla_x\Phi^f}(t,y)&:=\int_{\mR^{d_1}}\!c(x,y)\cdot\nabla_x\Phi^f(t,x,y)\mu^{y}(\dif x),\\
\overline{H\cdot\nabla_y\Phi^f}(t,y)&:=\int_{\mR^{d_1}}\!H(t,x,y)\cdot\nabla_y\Phi^f(t,x,y)\mu^{y}(\dif x).
\end{align*}
The following is the main result of this subsection, which will play a crucial role in the proof of Theorem \ref{main2}.

\bl\label{key22}
Let {\bf (A$_\sigma$)}, {\bf (A$_b$)}, (\ref{ac})  hold   and   $\delta\in(0,1]$. For any function $f$  satisfying (\ref{cen2}), we have:

\vspace{1mm}
\noindent(i) (Regime 1) if $b, \sigma\in C_b^{\delta,\vartheta}$ with $\vartheta\in(0,2]$, $c\in L^\infty_p$, $F, H, G\in L^\infty_p$ and $f\in C_p^{\vartheta/2,\delta,\vartheta}$, then
\begin{align*}
\mE\left(\frac{1}{\gamma_\eps}\int_0^tf(s,X_s^\eps,Y_s^\eps)\dif s\right)\leq C_t\Big(\frac{\alpha_\eps^{\vartheta}}{\gamma_\eps}+\frac{\alpha_\eps^2}{\gamma_\eps^2}+\frac{\alpha_\eps^2}{\beta_\eps\gamma_\eps}\Big);
\end{align*}
(ii) (Regime 2)  if $b, \sigma\in C_b^{\delta,\vartheta}$  with $\vartheta\in(0,2]$, $c\in C_p^{\delta,\vartheta}$, $F, H, G\in L^\infty_p$ and $f\in C_p^{\vartheta/2,\delta,\vartheta}$, then
\begin{align*}
\mE\!\left(\!\frac{1}{\gamma_\eps}\!\int_0^tf(s,X_s^\eps,Y_s^\eps)\dif s\!\!\right)\!-\mE\left(\int_0^t\overline{c\cdot\nabla_x\Phi^f}(s,Y_s^\eps)\dif s\right)\!\!\leq\!\! C_t\Big(\frac{\alpha_\eps^{\vartheta}}{\gamma_\eps}+\frac{\alpha_\eps^2}{\gamma_\eps^2}+\frac{\alpha_\eps^2}{\beta_\eps}\Big);
\end{align*}
(iii) (Regime 3) If $b, \sigma\in C_b^{\delta,1+\vartheta}$ with $\vartheta\in(0,1]$, $c\in L^\infty_p$, $F, G\in L^\infty_p $, $H\in C_{p}^{\vartheta/2,\delta,\vartheta}$ and $f\in C_p^{(1+\vartheta)/2,\delta,1+\vartheta}$, then
\begin{align*}
\mE\left(\frac{1}{\gamma_\eps}\int_0^tf(s,X_s^\eps,Y_s^\eps)\dif s\right)-\mE\left(\int_0^t\overline{H\cdot\nabla_y\Phi^f}(s,Y_s^\eps)\dif s\right)\leq C_t\Big(\alpha_\eps^{\vartheta}+\frac{\alpha_\eps}{\beta_\eps}\Big),
\end{align*}
(iv) (Regime 4) If $b, \sigma\in C_b^{\delta,1+\vartheta}$ with $\vartheta\in(0,1]$, $c\in C_p^{\delta,\vartheta}$, $F, G\in L^\infty_p$, $H\in C_{p}^{\vartheta/2,\delta,\vartheta}$ and $f\in C_p^{(1+\vartheta)/2,\delta,1+\vartheta}$, then
\begin{align*}
\mE\!\left(\!\frac{1}{\gamma_\eps}\int_0^t\!f(s,X_s^\eps,Y_s^\eps)\dif s\!\right)\!-\!\mE\left(\int_0^t\!\Big[\overline{c\cdot\nabla_x\Phi^f}+\overline{H\cdot\nabla_y\Phi^f}\Big](s,Y_s^\eps)\dif s\right)\!\leq\! C_t \alpha_\eps^{\vartheta},
\end{align*}
where $C_t>0$ is a constant independent on $\delta, \eps$.
\el

\begin{proof}
	As in the proof of Lemma \ref{key}, by (\ref{ito}) we have that
	\begin{align*}
	\hat\cQ(\eps)&:=\mE\left(\frac{1}{\gamma_\eps}\int_0^tf(s,X_s^\eps,Y_s^\eps)\dif s\right) =\frac{\alpha_\eps^2}{\gamma_\eps}\mE\big[\Phi_n^f(0,x,y)-\Phi_n^f({t},X_{t}^\eps,Y_{t}^\eps)\big]\\
	&\quad+\frac{\alpha_\eps^2}{\gamma_\eps}\mE\left(\int_0^{t}\big(\p_s+\sL_1\big)\Phi_n^f(s,X_s^\eps,Y_s^\eps)\dif s\right)\\
	&\quad+\frac{1}{\gamma_\eps}\mE\left(\int_0^{t}\big(\sL_0\Phi_n^f-\sL_0\Phi^f\big)(s,X_s^\eps,Y_s^\eps)\dif s\right)\\
	&\quad+\frac{\alpha_\eps^2}{\gamma_\eps^2}\mE\left(\int_0^{t}\sL_2\Phi_n^f(s,X_s^\eps,Y_s^\eps)\dif s\right)\\
	&\quad+\frac{\alpha_\eps^2}{\beta_\eps\gamma_\eps}\mE\left(\int_0^{t}\sL_3\Phi_n^f(s,X_s^\eps,Y_s^\eps)\dif s\right)=:\sum_{i=1}^5\hat\cQ_i(\eps),
	\end{align*}
	where $\Phi_n^f$ is the mollifyer of $\Phi^f$ defined as in (\ref{fn}).
	Below, we first prove the most general case {\it (iv)}, and then provide the proof of the other cases with sight  changes.
	
	\vspace{1mm}
	\noindent{\it (iv) (Regime 4)} In this case, according to Theorem \ref{popde}, (\ref{pro}) and by the assumptions on $b,\sigma$ and $f$, we have $\Phi^f\in C_p^{(1+\vartheta)/2,2+\delta,1+\vartheta}$ with $\vartheta\in (0,1]$. Following exactly the same arguments as in Lemma \ref{key}, we  get
	$$
	\hat\cQ_1(\eps)\leq C_1\frac{\alpha_\eps^2}{\gamma_\eps},
	$$
	and by (\ref{n222}),
	$$
	\hat\cQ_{2}(\eps)\leq C_2\frac{\alpha_\eps^2}{\gamma_\eps}n^{2-(1+\vartheta)}\mE\left(\int_0^t\big(1+|X_s^\eps|^{2m}\big)\dif s\right)\leq C_2\frac{\alpha_\eps^2}{\gamma_\eps}n^{1-\vartheta}.
	$$
	Meanwhile, we have by (\ref{n111}) that
	$$
	\hat\cQ_{3}(\eps)\leq C_3\frac{1}{\gamma_\eps}n^{-1-\vartheta}\mE\left(\int_0^t\big(1+|X_s^\eps|^{2m}\big)\dif s\right)\leq C_3\frac{1}{\gamma_\eps}n^{-1-\vartheta}.
	$$
	To control the last two terms, recall  that we have $\alpha_\eps=\beta_\eps=\gamma_\eps$ in this case. Thus, by the definition of $\sL_2$, we can write
	\begin{align*}
	&\hat\cQ_4(\eps)-\mE\left(\int_0^{t}\overline {H\cdot\nabla_y\Phi^f}(s,Y_s^\eps)\dif s\right)\\
	&\leq \mE\bigg(\!\int_0^{t}\!H(s,X_s^\eps,Y_s^\eps)\cdot\nabla_y\Phi_n^f(s,X_s^\eps,Y_s^\eps)\\
	&\quad-H(s,X_s^\eps,Y_s^\eps)\cdot\nabla_y\Phi^f(s,X_s^\eps,Y_s^\eps)\dif s\bigg)\\
	&\quad+ \mE\left(\int_0^{t}\!H(s,X_s^\eps,Y_s^\eps)\cdot\nabla_y\Phi^f(s,X_s^\eps,Y_s^\eps)-\overline {H\cdot\nabla_y\Phi^f}(s,Y_s^\eps)\dif s\right)\\
	&=:\hat\cQ_{41}(\eps)\!+\!\hat\cQ_{42}(\eps).
	\end{align*}	
	Since $\nabla_y\Phi^f\in C_p^{(1+\vartheta)/2,2+\delta,\vartheta}$ and by the fact
	$$
	\nabla_y(\Phi^f_n)=(\nabla_y\Phi^f)*\rho_1^n*\rho_2^n,
	$$
	we can deduce by using (\ref{n111}) again that
	\begin{align*}
	\hat\cQ_{41}(\eps)&\leq C_4\mE\left(\int_0^{t}\big\|\nabla_y\Phi_n^f(\cdot,X_s^\eps,\cdot)-\nabla_y\Phi^f(\cdot,X_s^\eps,\cdot)\big\|_\infty\big(1+|X_s^\eps|^m\big)\dif s\right)\\
	&\leq C_4n^{-\vartheta}\mE\left(\int_0^t\big(1+|X_s^\eps|^{2m}\big)\dif s\right)\leq C_4n^{-\vartheta}.
	\end{align*}	
	On the other hand, by the assumption that $H\in C_p^{\vartheta/2,\delta,\vartheta}$, it is easy to check that $H\cdot\nabla_y\Phi^f\in C_p^{\vartheta/2,\delta,\vartheta}$. Thus $\overline{H\cdot\nabla_y\Phi^f}\in C_b^{\vartheta/2,\vartheta}$ by Lemma \ref{aaav}. Note  that the function
	$$
	H(t,x,y)\cdot\nabla_y\Phi^f(t,x,y)-\overline {H\cdot\nabla_y\Phi^f}(t,y)
	$$
	satisfies (\ref{cen2}). Hence by applying Lemma \ref{key} with $\theta\in(0,1]$ it follows  that
	$$
	\hat\cQ_{42}(\eps)\leq C_4\Big(\alpha_\eps^{\vartheta}+\frac{\alpha_\eps^2}{\beta_\eps}\Big).
	$$
	Finally, due to the assumption that $c\in C_p^{\delta,\vartheta}$ and by the same idea as above, one can check that
	\begin{align*}
	\hat\cQ_5(\eps)-\mE\left(\int_0^{t}\overline {c\cdot\nabla_x\Phi^f}(s,Y_s^\eps)\dif s\right)\leq C_5\big(n^{-1-\vartheta}+\alpha_\eps^\vartheta\big).
	\end{align*}
	Combining the above computations, we arrive at
	\begin{align*} &\hat\cQ(\eps)-\mE\left(\int_0^t\!\Big[\overline{c\cdot\nabla_x\Phi^f}+\overline{H\cdot\nabla_y\Phi^f}\Big](s,Y_s^\eps)\dif s\right)\\
&\leq C_6\Big(\frac{\alpha_\eps^2}{\gamma_\eps}n^{1-\vartheta}+\frac{1}{\gamma_\eps}n^{-1-\vartheta}+n^{-\vartheta} +\alpha_\eps^{\vartheta}+\frac{\alpha_\eps^2}{\beta_\eps}\Big).
	\end{align*}
	Taking $n=\alpha_\eps^{-1}$, we get the desired result.
	
	\vspace{1mm}
	\noindent{\it (iii) (Regime 3)} The only difference to  {\it (iv)} is the last term. Recall that we have $\alpha_\eps=\gamma_\eps$ in this case. Thus, as above we have
	\begin{align*} &\hat\cQ_{1}(\eps)+\hat\cQ_{2}(\eps)+\hat\cQ_{3}(\eps)+\bigg[\hat\cQ_{4}(\eps)-\mE\left(\int_0^t\overline{H\cdot\nabla_y\Phi^f}(s,Y_s^\eps)\dif s\right)\bigg]\\
&\leq\! C_1\Big(\alpha_\eps n^{1-\vartheta}\!+\!\frac{1}{\alpha_\eps}n^{-1-\vartheta}\!+\!n^{-\vartheta} \!+\!\alpha_\eps^{\vartheta}\!+\!\frac{\alpha_\eps^2}{\beta_\eps}\Big).
	\end{align*}
	To control $\hat\cQ_5(\eps)$, we can simply use the growth condition on $c$ to get
	$$
	\hat\cQ_5(\eps)\leq C_2\frac{\alpha_\eps^2}{\beta_\eps\gamma_\eps} \mE\left(\int_0^t\big(1+|X_s^\eps|^m\big)\dif s\right)\leq C_2\frac{\alpha_\eps}{\beta_\eps}.
	$$
	Consequently, we arrive at
	\begin{align*}
	&\hat\cQ(\eps)-\mE\left(\int_0^t\overline{H\cdot\nabla_y\Phi^f}(s,Y_s^\eps)\dif s\right)\\
&\leq C_3\Big(\frac{1}{\alpha_\eps}n^{-1-\vartheta}+\alpha_\eps n^{1-\vartheta}+n^{-\vartheta} +\alpha_\eps^{\vartheta}+\frac{\alpha_\eps}{\beta_\eps}\Big)
	\end{align*}
	Taking $n=\alpha_\eps^{-1}$ again, we get the desired result.
	
	\vspace{1mm}
	\noindent{\it (ii) (Regime 2)} Now we need to pay attention to the assumption that $\vartheta\in (0,2]$. By Theorem \ref{popde}, (\ref{pro}) and the assumptions on $b,\sigma$ and $f$, we have $\Phi^f\in C_p^{\vartheta/2,2+\delta,\vartheta}$. Arguing as before, we   have
	$$
	\hat\cQ_1(\eps)+\hat\cQ_2(\eps)+\hat\cQ_3(\eps)\leq C_1\Big(\frac{\alpha_\eps^2}{\gamma_\eps}+\frac{\alpha_\eps^2}{\gamma_\eps}n^{2-\vartheta}+\frac{1}{\gamma_\eps}n^{-\vartheta}\Big),
	$$
	and for $\hat\cQ_4(\eps)$, as in the proof of Lemma \ref{key}, we have
	$$
	\hat\cQ_4(\eps)\leq C_2\frac{\alpha_\eps^2}{\gamma_\eps^2} n^{1-(\vartheta\wedge1)}\mE\left(\int_0^t\big(1+|X_s^\eps|^{2m}\big)\dif s\right)\leq C_2\frac{\alpha_\eps^2}{\gamma_\eps^2} n^{1-(\vartheta\wedge1)}.
	$$
	For the last term, recall that $\alpha_\eps^2=\beta_\eps\gamma_\eps$ in this case. We   use Lemma \ref{key} with $\vartheta\in(0,2]$ to deduce that
	\begin{align*}
	&\hat\cQ_5(\eps)-\mE\left(\int_0^{t}\overline {c\cdot\nabla_x\Phi^f}(s,Y_s^\eps)\dif s\right)\\
	&\leq \mE\left(\int_0^{t}c(X_s^\eps,Y_s^\eps)\cdot\big[\nabla_x\Phi_n^f(s,X_s^\eps,Y_s^\eps)-\nabla_x\Phi^f(s,X_s^\eps,Y_s^\eps)\big]\dif s\right)\\
	&\quad+\mE\left(\int_0^{t}c(X_s^\eps,Y_s^\eps)\cdot\nabla_x\Phi^f(s,X_s^\eps,Y_s^\eps)-\overline {c\cdot\nabla_x\Phi^f}(s,Y_s^\eps)\dif s\right)\\ &\leq C_3\Big(n^{-\vartheta}+\alpha_\eps^{\vartheta}+\alpha_\eps^{\vartheta\wedge1}\cdot\frac{\alpha_\eps}{\gamma_\eps}+\frac{\alpha_\eps^2}{\beta_\eps}\Big).
	\end{align*}
	Consequently, we arrive at
	\begin{align*}
	&\hat\cQ(\eps)-\mE\left(\int_0^{t}\overline {c\cdot\nabla_x\Phi^f}(s,Y_s^\eps)\dif s\right)\\
&\leq C_4\Big(\frac{\alpha_\eps^2}{\gamma_\eps}n^{2-\vartheta}+\frac{1}{\gamma_\eps}n^{-\vartheta}+\frac{\alpha_\eps^2}{\gamma_\eps^2} n^{1-(\vartheta\wedge1)}+\alpha_\eps^{\vartheta}+\alpha_\eps^{\vartheta\wedge1}\cdot\frac{\alpha_\eps}{\gamma_\eps}+\frac{\alpha_\eps^2}{\beta_\eps}\Big).
	\end{align*}
	Still, we can take $n=\alpha_\eps^{-1}$ to get that for $\vartheta\in(0,1]$,
	\begin{align*}
	&\hat\cQ(\eps)-\mE\left(\int_0^{t}\overline {c\cdot\nabla_x\Phi^f}(s,Y_s^\eps)\dif s\right)\\
&\leq C_4\Big(\frac{\alpha_\eps^\vartheta}{\gamma_\eps}+\frac{\alpha_\eps^{1+\vartheta}}{\gamma_\eps^2} +\frac{\alpha_\eps^2}{\beta_\eps}\Big)\leq C_4\Big(\frac{\alpha_\eps^\vartheta}{\gamma_\eps} +\frac{\alpha_\eps^2}{\beta_\eps}\Big).
	\end{align*}
	and for $\vartheta\in(1,2]$, we have
	\begin{align*}
	\hat\cQ(\eps)-\mE\left(\int_0^{t}\overline {c\cdot\nabla_x\Phi^f}(s,Y_s^\eps)\dif s\right)\leq C_4\Big(\frac{\alpha_\eps^\vartheta}{\gamma_\eps}+\frac{\alpha_\eps^{2}}{\gamma_\eps^2} +\frac{\alpha_\eps^2}{\beta_\eps}\Big),
	\end{align*}
	which in turn yields the desired result.
	
	\vspace{1mm}
	\noindent
	{\it (i) (Regime 1)} Finally, we can use Lemma \ref{key} directly to get that
	\begin{align*}
	\mE\left(\frac{1}{\gamma_\eps}\int_0^tf(s,X_s^\eps,Y_s^\eps)\dif s\right)&\leq C_1\Big(\frac{\alpha_\eps^{\vartheta}}{\gamma_\eps}+\frac{\alpha_\eps^{1+(\vartheta\wedge1)}}{\gamma_\eps^2}+\frac{\alpha_\eps^2}{\beta_\eps\gamma_\eps}\Big)\\
	&\leq C_1\Big(\frac{\alpha_\eps^{\vartheta}}{\gamma_\eps}+\frac{\alpha_\eps^2}{\gamma_\eps^2}+\frac{\alpha_\eps^2}{\beta_\eps\gamma_\eps}\Big).
	\end{align*}
	The whole proof is finished.
\end{proof}

\section{Diffusion approximations}
Recall that $\Phi(t,x,y)$ is the solution to Poisson equation (\ref{pde1}), and the limit effective system is given by (\ref{sde000}).
For $k=1,\cdots,4$, denote by $\hat \cL_k$ the infinitesimal operator of $\hat Y^k_t$, i.e.,
\begin{align*}
\hat\cL_k:=\sum_{i,j=1}^{d_2} \hat \cG^{ij}_k(t,y)\frac{\p^2}{\p y_i\p y_j}+\sum_{i=1}^{d_2}\hat F^i_k(t,y)\frac{\p}{\p y_i},
\end{align*}
where $\hat \cG_k(t,y):=\hat G_k\hat G^*_k(t,y)/2$. The following  properties for the averaged coefficients follow  by Lemma \ref{aaav}.

\bl\label{ar}
Under the assumptions in Theorem \ref{main2}, for every $k=1,\cdots,4$, $\hat\cG_k$ is non-degenerate in $y$ uniformly with respect to $t$. Moreover, we have $\hat F_k, \hat \cG_k\in C_b^{\vartheta/2,\vartheta}$.
\el

\begin{proof}
	The non-degeneracy of $\hat\cG_1=\hat\cG_2$ follows directly by {\bf (A$_G$)} and the definition.
	Furthermore, by (\ref{pro})  we have
	\begin{align*}
	&\int_{\mR^{d_1}}H(t,x,y)\Phi^*(t,x,y)\mu^y(\dif x)\\
	&=\mE\left(\int_0^\infty\int_{\mR^{d_1}}H(t,x,y) H^*(t,X_t^y,y)\mu^y(\dif x)\dif t\right),
	\end{align*}
	which is non-negative   by the homogeneity of $X_t^y$. Thus $\hat \cG_3=\hat\cG_4$ are also non-degenerate. Let us prove the regularity for the most general case $\hat F_4$. The conclusions  for $\hat F_2, \hat F_3, \hat F_4$ and $\hat G_k$ ($k=1,\cdots,4$) can be proved by the same argument.
	In this case, we have by the assumption that $b,\sigma\in C_b^{\delta,1+\vartheta}$, $H\in C_p^{(1+\vartheta)/2,\delta,1+\vartheta}$ and Theorem \ref{popde} that $\Phi(t,x,y)\in C_p^{(1+\vartheta)/2,2+\delta,1+\vartheta}$, which together with the assumptions on  $F$ and $c$ implies that
	$$
	F(t,x,y)+c(x,y)\cdot\nabla_x\Phi(t,x,y)+H(t,x,y)\cdot\nabla_y\Phi(t,x,y)\in C_p^{\vartheta/2,\delta,\vartheta}.
	$$
	Thus the conclusion follows by Lemma \ref{aaav} directly.
\end{proof}

According to the above result, there exists a unique weak solution $\hat Y_t^k$ for SDE (\ref{sde000}) for every $k=1,\cdots,4$. To prove the weak convergence of $Y_t^\eps$ to $\hat Y_t^k$, we need to consider the following Cauchy problem on $[0,T]\times\mR^{d_2}$:
\begin{equation}\left\{\begin{array}{l}\label{PDE}
\displaystyle
\p_t u_k(t,y)-\hat\cL_k  u_k(t,y)=0,\quad t\in [0, T),\\
u_k(0,y)=\varphi(y),
\end{array}\right.
\end{equation}
where $T>0$ and  $\varphi$ is a function on $\mR^{d_2}$.
As a direct consequence of Lemma \ref{ar}, we have the following result for the Cauchy problem (\ref{PDE}), which is well-known in the theory of PDEs, see e.g. \cite[Chapter IV, Section 5]{La-So-Ur}.

\bl\label{ppp}
Assume that $\varphi\in C_b^{2+\vartheta}$. Then for every $k=1,\cdots, 4$, there exists a unique solution $u_k\in C_b^{(2+\vartheta)/2,2+\vartheta}$ to equation (\ref{PDE}) which is given by
\begin{align}\label{cau}
u_k(t,y)=\mE\varphi(\hat Y^k_t(y)).
\end{align}
Moreover, we also have $\nabla_yu_k\in C_b^{(1+\vartheta)/2,1+\vartheta}$ and $\nabla^2_yu\in C_b^{\vartheta/2,\vartheta}$.
\el

Now, we are in the position to give:

\begin{proof}[Proof of Theorem \ref{main2}]
	Given $T>0$ and  $\varphi\in C_b^{2+\vartheta}$, let $u_k$ ($k=1,\cdots, 4$) be the unique solution to equation (\ref{PDE}). Define
	$$
	\tilde u_k(t,y):=u_k(T-t,y).
	$$
	Then it is obvious that $\tilde u_k(T,y)=u_k(0,y)=\varphi(y)$. As a result, we can deduce by (\ref{cau}) and the It\^o's formula that
	\begin{align*}
	&\cR_k(\eps):=\mE[\varphi(Y_T^\eps)]-\mE[\varphi(\hat Y^k_T)]=\mE[\tilde u_k(T,Y_T^\eps)-\tilde u_k(0,y)]\\
	&=\mE\!\left(\int_0^T\big(\p_s+\sL_1\big)\tilde u_k(s,Y_s^\eps)\dif s\right)+\mE\left(\frac{1}{\gamma_\eps}\int_0^T\!\sL_2\tilde u_k(s,Y_s^\eps)\dif s\right)\\
	&=\mE\!\left(\int_0^T\!\!\big(\sL_1-\hat\cL_k\big)\tilde u_k(s,Y_s^\eps)\dif s\right)\!+\!\frac{1}{\gamma_\eps}\mE\!\left(\int_0^T\!\! H(s,X_s^\eps,Y_s^\eps)\nabla_y\tilde u_k(s,Y_s^\eps)\dif s\right),
	\end{align*}
	where $\sL_2$ and $\sL_1$ are defined by (\ref{lll}). Below, we divide the proof into four parts, which correspond  to the four regimes in  (\ref{regime}).
	
	\vspace{1mm}\noindent
	(i) (Regime 1) In this case, we have
	\begin{align*}
	\big(\sL_1-\hat\cL_1\big)\tilde u_1(t,y)&=\big[\cG(t,x,y)-\hat\cG_1(t,y)\big]\cdot\nabla^2_y\tilde u_1(t,y)\\
	&\quad+\big[F(t,x,y)-\hat F_1(t,y)\big]\cdot\nabla_y\tilde u_1(t,y).
	\end{align*}
	Note that the function $\big(\sL_1-\hat\cL_1\big)\tilde u_1(t,y)$ satisfies the centering condition (\ref{cen2}). Moreover, by the assumption that $F, \cG\in C_p^{\vartheta/2,\delta,\vartheta}$, Lemma \ref{ar} and Lemma \ref{ppp}, we have $\big(\sL_1-\hat\cL_1\big)\tilde u_1(t,y) \in C_p^{\vartheta/2,\delta,\vartheta}$.
	Thus,
	applying  Lemma \ref{key} with $\vartheta\in(0,2]$ we can get
	\begin{align*}
	\cR_{11}(\eps):=\mE\left(\int_0^T\big(\sL_1-\hat\cL_1\big)\tilde u_1(s,Y_s^\eps)\dif s\right)\leq C_1\Big(\alpha_\eps^{\vartheta}+\alpha_\eps^{\vartheta\wedge1}\cdot\frac{\alpha_\eps}{\gamma_\eps}+\frac{\alpha_\eps^2}{\beta_\eps}\Big).
	\end{align*}
	On the other hand, note that $H(t,x,y)\cdot\nabla_y\tilde u_1(t,y)\in C_p^{\vartheta/2,\delta,\vartheta}$ also satisfies (\ref{cen2}). Hence we   have by Lemma \ref{key22} {\it (i)} that
	\begin{align*}
	\cR_{12}(\eps)&:=\frac{1}{\gamma_\eps}\mE\left(\int_0^T H(s,X_s^\eps,Y_s^\eps)\cdot\nabla_y\tilde u_1(s,Y_s^\eps)\dif s\right)\\
	&\leq C_1\Big(\frac{\alpha_\eps^{\vartheta}}{\gamma_\eps}+\frac{\alpha_\eps^2}{\gamma_\eps^2}+\frac{\alpha_\eps^2}{\beta_\eps\gamma_\eps}\Big).
	\end{align*}
	Consequently,
	$$
	\cR_1(\eps)=\cR_{11}(\eps)+\cR_{12}(\eps)\leq C_1\Big(\frac{\alpha_\eps^{\vartheta}}{\gamma_\eps}+\frac{\alpha_\eps^2}{\gamma_\eps^2}+\frac{\alpha_\eps^2}{\beta_\eps\gamma_\eps}\Big).
	$$
	(ii) (Regime 2) In this case, note that
	\begin{align*}
	\hat\cL_2=\hat\cL_1+\overline{c\cdot\nabla_x\Phi}(t,y)\cdot\nabla_y.
	\end{align*}
	Thus, we can write
	\begin{align*}
	\cR_2(\eps)&=\mE\left(\int_0^T\big(\sL_1-\hat\cL_1\big)\tilde u_2(s,Y_s^\eps)\dif s\right)\\
	&\quad+\bigg[\frac{1}{\gamma_\eps}\mE\left(\int_0^T H(s,X_s^\eps,Y_s^\eps)\cdot\nabla_y\tilde u_2(s,Y_s^\eps)\dif s\right)\\
	&\quad-\mE\left(\int_0^T\overline{c\cdot\nabla_x\Phi}(s,Y_s^\eps)\cdot\nabla_y\tilde u_2(s,Y_s^\eps)\dif s\right)\bigg]=:\cR_{21}(\eps)+\cR_{22}(\eps).
	\end{align*}
	Entirely similar as above, we can control the first term by
	\begin{align*}
	\cR_{21}(\eps)\leq C_2\Big(\alpha_\eps^{\vartheta}+\alpha_\eps^{\vartheta\wedge1}\cdot\frac{\alpha_\eps}{\gamma_\eps}+\frac{\alpha_\eps^2}{\beta_\eps}\Big).
	\end{align*}
	On the other hand, let $\Psi_2$ be the solution to
	$$
	\sL_0(x,y)\Psi_2(t,x,y)=-H(t,x,y)\cdot\nabla_y\tilde u_2(t,y)\in C_p^{\vartheta/2,\delta,\vartheta}.
	$$
	Then, it is easy to see that
	$$
	\Psi_2(t,x,y)=\Phi(t,x,y)\cdot\nabla_y\tilde u_2(t,y),
	$$
	where $\Phi$ satisfies (\ref{pde1}).
	Furthermore,
	$$
	\int_{\mR^{d_1}}c(x,y)\cdot\nabla_x\Psi_2(t,x,y)\mu^y(\dif x)=\overline{c\cdot\nabla_x\Phi}(t,y)\cdot\nabla_y\tilde u_2(t,y).
	$$
	Consequently, we can apply Lemma \ref{key22} {\it (ii)} with $f=H\cdot\nabla_y\tilde u_2$ to get that
	\begin{align*}
	\cR_{22}(\eps)\leq C_2\Big(\frac{\alpha_\eps^{\vartheta}}{\gamma_\eps}+\frac{\alpha_\eps^2}{\gamma_\eps^2}+\frac{\alpha_\eps^2}{\beta_\eps}\Big).
	\end{align*}
	Thus we arrive at
	\begin{align*}
	\cR_{2}(\eps)=\cR_{21}(\eps)+\cR_{22}(\eps)\leq C_2\Big(\frac{\alpha_\eps^{\vartheta}}{\gamma_\eps}+\frac{\alpha_\eps^2}{\gamma_\eps^2}+\frac{\alpha_\eps^2}{\beta_\eps}\Big).
	\end{align*}
	(iii) (Regime 3) In this case, we have
	\begin{align*}
	\hat\cL_3=\hat\cL_1+\overline{H\cdot\nabla_y\Phi}(t,y)\cdot\nabla_y+\overline{H\cdot\Phi^*}(t,y)\cdot\nabla^2_y,
	\end{align*}
	where the extra diffusion coefficient is given by
	\begin{align*}
	\overline{H\cdot\Phi^*}(t,y):=\int_{\mR^{d_1}}H(t,x,y)\Phi^*(t,x,y)\mu^y(\dif x).
	\end{align*}
	Thus we can write
	\begin{align*}
	&\cR_3(\eps)=\mE\left(\int_0^T\big(\sL_1-\hat\cL_1\big)\tilde u_3(s,Y_s^\eps)\dif s\right)\\
	&+\bigg[\frac{1}{\gamma_\eps}\mE\left(\int_0^T H(s,X_s^\eps,Y_s^\eps)\cdot\nabla_y\tilde u_3(s,Y_s^\eps)\dif s\right)\\
	&-\mE\left(\int_0^T\overline{H\cdot\nabla_y\Phi}(s,Y_s^\eps)\nabla_y\tilde u_3(s,Y_s^\eps)+\overline{H\cdot\Phi^*}(s,Y_s^\eps)\nabla^2_y\tilde u_3(s,Y_s^\eps)\dif s\right)\bigg]\\
	&=:\cR_{31}(\eps)+\cR_{32}(\eps).
	\end{align*}
	Following the same idea as in (i) and using Lemma \ref{key} with $\vartheta\in(0,1]$, we can control the first term by
	\begin{align*}
	\cR_{31}(\eps)\leq C_3\Big(\alpha_\eps^{\vartheta}+\frac{\alpha_\eps^2}{\beta_\eps}\Big).
	\end{align*}
	On the other hand, let $\Psi_3$ be the solution to
	$$
	\sL_0(x,y)\Psi_3(t,x,y)=-H(t,x,y)\cdot\nabla_y\tilde u_3(t,y)\in C_p^{(1+\vartheta)/2,\delta,1+\vartheta},
	$$
	which is given by
	$$
	\Psi_3(t,x,y)=\Phi(t,x,y)\cdot\nabla_y\tilde u_3(t,y).
	$$
	One can  check that
	\begin{align*}
	&\int_{\mR^{d_1}}H(t,x,y)\cdot\nabla_y\Psi_3(t,x,y)\mu^y(\dif x)\\
	&=\int_{\mR^{d_1}}H(t,x,y)\cdot\nabla_y\Phi(t,x,y)\mu^y(\dif x)\cdot\nabla_y\tilde u_3(t,y)\\
	&\quad+\int_{\mR^{d_1}}H(t,x,y)\cdot\Phi^*(t,x,y)\mu^y(\dif x)\cdot\nabla^2_y\tilde u_3(t,y)\\
	&=\overline{H\cdot\nabla_y\Phi}(t,y)\cdot\nabla_y\tilde u_3(t,y)+\overline{H\cdot\Phi^*}(t,y)\cdot\nabla^2_y \tilde u_3(t,y).
	\end{align*}
	Consequently, by applying Lemma \ref{key22} {\it (iii)} with $f=H\cdot\nabla_y\tilde u_3$ we have that
	\begin{align*}
	\cR_{32}(\eps)\leq C_3\Big(\alpha_\eps^{\vartheta}+\frac{\alpha_\eps}{\beta_\eps}\Big),
	\end{align*}
	which in turn yields the desired result.
	
	\vspace{1mm}\noindent
	(iv) (Regime 4) Finally, we have
	\begin{align*}
	\hat\cL_4=\hat\cL_1+\overline{c\cdot\nabla_x\Phi}(t,y)\cdot\nabla_y+\overline{H\cdot\nabla_y\Phi}(t,y)\cdot\nabla_y +\overline{H\cdot\Phi^*}(t,y)\cdot\nabla^2_y.
	\end{align*}
	We thus have
	\begin{align*}
	\cR_4(\eps)&=\mE\left(\int_0^T\big(\sL_1-\hat\cL_1\big)\tilde u_4(s,Y_s^\eps)\dif s\right)\\
	&\quad+\bigg[\frac{1}{\gamma_\eps}\mE\left(\int_0^T H(s,X_s^\eps,Y_s^\eps)\cdot\nabla_y\tilde u_4(s,Y_s^\eps)\dif s\right)\\
	&\quad-\mE\bigg(\int_0^T\Big[\overline{c\cdot\nabla_x\Phi}+\overline{H\cdot\nabla_y\Phi}\Big](s,Y_s^\eps)\cdot\nabla_y\tilde u_4(s,Y_s^\eps)\\
	&\qquad\quad\quad+\overline{H\cdot\Phi^*}(s,Y_s^\eps)\cdot\nabla^2_y\tilde u_4(s,Y_s^\eps)\dif s\bigg)\bigg]=:\cR_{41}(\eps)+\cR_{42}(\eps).
	\end{align*}
	As above, we can control the first term by
	\begin{align*}
	\cR_{41}(\eps)\leq C_4\alpha_\eps^{\vartheta}.
	\end{align*}
	Let $\Psi_4$ be the solution to
	$$
	\sL_0(x,y)\Psi_4(t,x,y)=-H(t,x,y)\cdot\nabla_y\tilde u_4(t,y)\in C_p^{(1+\vartheta)/2,\delta,1+\vartheta},
	$$
	which is given by
	$$
	\Psi_4(t,x,y)=\Phi(t,x,y)\cdot\nabla_y\tilde u_4(t,y).
	$$
	One can  check that
	\begin{align*}
	&\int_{\mR^{d_1}}c(x,y)\cdot\nabla_x\Psi_4(t,x,y)\mu^y(\dif x)+\int_{\mR^{d_1}}H(t,x,y)\cdot\nabla_y\Psi_4(t,x,y)\mu^y(\dif x)\\
	&=\Big[\overline{c\cdot\nabla_x\Phi}+\overline{H\cdot\nabla_y\Phi}\Big](t,y)\cdot\nabla_y\tilde u_4(t,y)+\overline{H\cdot\Phi^*}(t,y)\cdot\nabla^2_y \tilde u_4(t,y).
	\end{align*}
	Consequently, by applying Lemma \ref{key22} {\it (iv)} with $f=H\cdot\nabla_y\tilde u_4$ we have that
	\begin{align*}
	\cR_{42}(\eps)\leq C_4\alpha_\eps^{\vartheta}.
	\end{align*}
	The whole proof is finished.
\end{proof}

\end{document}